\documentclass[a4paper]{article}

\usepackage[english]{babel}
\usepackage[utf8]{inputenc}
\usepackage{amsmath}
\usepackage{graphicx}
\usepackage[colorinlistoftodos]{todonotes}
\usepackage[]{natbib}
\usepackage{amssymb}
\usepackage{amsmath}
\usepackage{MnSymbol}
\usepackage{graphicx}
\usepackage{booktabs}
\usepackage[margin=20pt]{subcaption}
\usepackage{setspace}
\usepackage{blindtext}
\usepackage[utf8]{inputenc}
\usepackage{tabularx, ragged2e}
\usepackage{bbm}
\usepackage{latexsym}
\usepackage{caption}
\usepackage{float}
\usepackage{comment}
\usepackage{dcolumn}
\usepackage{adjustbox}
\usepackage{amsthm}
\usepackage{authblk}
\usepackage[colorlinks,citecolor=blue,colorlinks=true,urlcolor=blue]{hyperref}
\usepackage[normalem]{ulem}
\usepackage{dsfont}

\newtheorem{theorem}{Theorem}[section]
\newtheorem{corollary}{Corollary}[section]
\newtheorem{lemma}{Lemma}[section]

\newcommand{\cov}[1]{\text{Cov}{\left(#1\right)}}
\newcommand{\var}[1]{\text{Var}{\left(#1\right)}}
\newcommand{\diag}[1]{\textbf{Diag}{\left(#1\right)}}
\newcommand{\E}[1]{\mathbb{E}\left[#1\right]}

\title{Ridge Rerandomization: An Experimental Design Strategy in the Presence of Collinearity}

\author[1]{Zach Branson}

\author[2]{Stephane Shao\thanks{
    We would like to thank Espen Bernton and Tirthankar Dasgupta for insightful comments throughout the progress of this work. We would also like to especially thank Evan Greif for discussions that initially motivated this work. This research was supported by the National Science Foundation Graduate Research Fellowship Program under Grant No. 1144152. Any opinions, findings, and conclusions or recommendations expressed in this material are those of the authors and do not necessarily reflect the views of the National Science Foundation.}}

\affil[1]{Department of Statistics and Data Science, Carnegie Mellon University}
\affil[2]{Department of Statistics, Harvard University}

\date{\today}

\begin{document}
\maketitle

\allowdisplaybreaks
\begin{abstract}
	Randomization ensures that observed and unobserved covariates are balanced, on average. However, randomizing units to treatment and control often leads to covariate imbalances in realization, and such imbalances can inflate the variance of estimators of the treatment effect. One solution to this problem is rerandomization---an experimental design strategy that randomizes units until some balance criterion is fulfilled---which yields more precise estimators of the treatment effect if covariates are correlated with the outcome. Most rerandomization schemes in the literature utilize the Mahalanobis distance, which may not be preferable when covariates are correlated or vary in importance. As an alternative, we introduce an experimental design strategy called ridge rerandomization, which utilizes a modified Mahalanobis distance that addresses collinearities among covariates and automatically places a hierarchy of importance on the covariates according to their eigenstructure. This modified Mahalanobis distance has connections to principal components and the Euclidean distance, and---to our knowledge---has remained unexplored. We establish several theoretical properties of this modified Mahalanobis distance and our ridge rerandomization scheme. These results guarantee that ridge rerandomization is preferable over randomization and suggest when ridge rerandomization is preferable over standard rerandomization schemes. We also provide simulation evidence that suggests that ridge rerandomization is particularly preferable over typical rerandomization schemes in high-dimensional or high-collinearity settings.
\end{abstract}

\section{Introduction}
\label{sec:introduction}

Randomized experiments are often considered the ``gold standard" of scientific investigations because, on average, randomization balances all potential confounders, both observed and unobserved \citep{krause2003random}. However, many have noted that randomized experiments can yield ``bad allocations,'' where some covariates are not well-balanced across treatment groups \citep{seidenfeld1981levi,lindley1982role,papineau1994virtues,rosenberger2008handling}. Covariate imbalance among different treatment groups complicates the interpretation of estimated causal effects, and thus covariate adjustments are often employed, typically through regression or other comparable methods.

However, it would be better to prevent such covariate imbalances from occurring before treatment is administered, rather than depend on assumptions for covariate adjustment post-treatment which may not hold \citep{freedman2008regression}. One common experimental design tool is blocking, where units are first grouped together based on categorical covariates, and then treatment is randomized within these groups. However, blocking is less intuitive when there are non-categorical covariates. A more recent experimental design tool that prevents covariate imbalance and allows for non-categorical covariates is the rerandomization scheme of \cite{morgan2012rerandomization}, where units are randomized until a prespecified level of covariate balance is achieved. Rerandomization has been discussed as early as R.A. Fisher (e.g., see \cite{fisher1992arrangement}), and more recent works (e.g., \cite{cox2009randomization}, \cite{bruhn2009pursuit}, and \cite{worrall2010evidence}) recommend rerandomization. \cite{morgan2012rerandomization} formalized these recommendations in treatment-versus-control settings and was one of the first works to establish a theoretical framework for rerandomization schemes. Since \cite{morgan2012rerandomization}, several extensions have been made. \cite{morgan2015rerandomization} developed rerandomization for treatment-versus-control experiments where there are tiers of covariates that vary in importance; \cite{branson2016improving} extended rerandomization to $2^K$ factorial designs; and \cite{ernst2017sequential} developed a rerandomization scheme for sequential designs. Finally, \cite{li2018asymptotic} established asymptotic results for the rerandomization schemes considered in \cite{morgan2012rerandomization} and \cite{morgan2015rerandomization}, and \cite{li2020rerandomization} established asymptotic results for regression adjustment combined with rerandomization.

All of these works focus on using an omnibus measure of covariate balance---the Mahalanobis distance \citep{mahalanobis1936generalised}---during the rerandomization scheme. The Mahalanobis distance is well-known within the matching and observational study literature, where it is used to find subsets of the treatment and control that are similar \citep{rubin1974multivariate,rosenbaum1985constructing,gu1993comparison,rubin2000combining}. The Mahalanobis distance is particularly useful in rerandomization schemes because (1) it is symmetric in the treatment assignment, which leads to unbiased estimators of the average treatment effect under rerandomization; and (2) it is equal-percent variance reducing if the covariates are ellipsoidally symmetric, meaning that rerandomization using the Mahalanobis distance reduces the variance of all covariate mean differences by the same percentage \citep{morgan2012rerandomization}.

However, the Mahalanobis distance is known to perform poorly in matching for observational studies when covariates are not ellipsoidally symmetric, there are strong collinearities among the covariates, or there are many covariates \citep{gu1993comparison,olsen1997multivariate,stuart2010matching}. One reason for this is that matching using the Mahalanobis distance places equal importance on balancing all covariates as well as their interactions \citep{stuart2010matching}, and this issue also occurs in rerandomization schemes that use the Mahalanobis distance. This issue was partially addressed by \cite{morgan2015rerandomization}, who proposed an extension of \cite{morgan2012rerandomization} that incorporates tiers of covariates that vary in importance, such that the most important covariates receive the most variance reduction. However, this requires researchers to specify an explicit hierarchy of importance for the covariates, which might be difficult, especially when the number of covariates is large. Furthermore, it is unclear how to conduct current rerandomization schemes if collinearity is so severe that the covariance matrix of covariates is degenerate, and thus the Mahalanobis distance is undefined.

As an alternative, we consider a rerandomization scheme using a modified Mahalanobis distance that inflates the eigenvalues of the covariates' covariance matrix to alleviate collinearities among the covariates, which has connections to ridge regression \citep{hoerl1970ridge}. Such a quantity has remained largely unexplored in the literature. First we establish several theoretical properties about this quantity, as well as several properties about a rerandomization scheme that uses this quantity. We show through simulation that a rerandomization scheme that incorporates this modified criterion can be beneficial in terms of variance reduction when there are strong collinearities among the covariates. In particular, this rerandomization scheme automatically specifies a hierarchy of importance based on the eigenstructure of the covariates, which can be useful when researchers are unsure about how much importance they should place on each covariate when designing a randomized experiment. We also discuss how this modified Mahalanobis distance connects to other criteria, such as principal components and the Euclidean distance. Because the rerandomization literature has focused almost exclusively on the Mahalanobis distance, this work also contributes to the literature by exploring the use of other criteria besides the Mahalanobis distance for rerandomization schemes.

The remainder of this paper is organized as follows. In Section \ref{sec:notation}, we introduce the notation that will be used throughout the paper. In Section \ref{sec: review}, we review the rerandomization scheme of \cite{morgan2012rerandomization}. In Section \ref{s:ridgeRerandomization}, we outline our proposed rerandomization approach and establish several theoretical properties of this approach, as well as several theoretical properties about the modified Mahalanobis distance. In Section \ref{s:simulations}, we provide simulation evidence that suggests that our rerandomization approach is often preferable over other rerandomization approaches, particularly in high-dimensional or high-collinearity settings. In Section \ref{s:discussion}, we conclude with a discussion of future work.

\section{Notation}
\label{sec:notation}
We use the colon notation $\lambda_{1:K} = (\lambda_1,...,\lambda_K)\in\mathbb{R}^K$ for tuples of objects, and we let $f(\lambda_{1:K}) = (f(\lambda_1),...,f(\lambda_K))$ for any univariate function $f:\mathbb{R}\rightarrow\mathbb{R}$. We respectively denote by $\mathbf{I}_N$ and $\mathbf{1}_N$ the $N\times N$ identity matrix and the $N$-dimensional column vector whose coefficients are all equal to $1$. Given a matrix $A$, we denote by $A_{ij}$ its $(i,j)$-coefficient, $A_{i\bullet}$ its $i$-th row, $A_{\bullet j}$ its $j$-th column, $A^\top$ its transpose, and $\text{tr}(A)$ its trace when $A$ is square. Given two symmetric matrices $A$ and $B$ of the same size, we write $A>B$ (resp.\! $A\geq B$) if the matrix $A-B$ is positive definite (resp.\! semi-definite).

Let $\mathbf{x}$ be the $N\times K$ matrix representing $K$ covariates measured on $N$ experimental units. Let $W_i = 1$ if unit $i$ is assigned to treatment and 0 otherwise, and let $\mathbf{W} = (W_1\,\dots\, W_N)^\top$. Unless stated otherwise, we will focus on completely randomized experiments \citep[see Definition 4.2]{imbens2015causal} with a fixed number of $N_T$ treated units and $N_C=N-N_T$ control units. For a given assignment vector $\mathbf{W}$, we define $\bar{\mathbf{x}}_T = N_T^{-1\,} \mathbf{x}^\top \mathbf{W}$ and $\bar{\mathbf{x}}_C = N_C^{-1\,}\mathbf{x}^\top\left(\mathbf{1}_N-\mathbf{W}\right)$ as the respective covariate mean vectors within treatment and control.

For completely randomized experiments, the covariance matrix of the covariate mean differences is $\mathbf{\Sigma} = \cov{\bar{\mathbf{x}}_T-\bar{\mathbf{x}}_C\,|\,\mathbf{x}} = N N_T^{-1}N_C^{-1} S^2_{\mathbf{x}}$, where $S^2_{\mathbf{x}} = (N-1)^{-1}(\mathbf{x}-\mathbf{1}_N\bar{\mathbf{x}}_N)^\top(\mathbf{x}-\mathbf{1}_N\bar{\mathbf{x}}_N)$ is the sample covariance matrix of $\mathbf{x}$ with $\bar{\mathbf{x}}_N = N^{-1} \mathbf{1}_N^\top \mathbf{x}$ \citep{morgan2012rerandomization}. Throughout, we use $\mathbf{\Sigma}$ to refer to this fixed covariance matrix, and we assume $\mathbf{\Sigma}>0$. The spectral decomposition ensures that $\mathbf{\Sigma}$ is diagonalizable with eigenvalues $\lambda_1\geq ... \geq\lambda_K>0$. Let $\mathbf{\Gamma}$ be the orthogonal matrix of corresponding eigenvectors, so that we may write $\mathbf{\Sigma} = \mathbf{\Gamma}\diag{\lambda_{1:K}}\mathbf{\Gamma}^\top$, where $\diag{\lambda_{1:K}}$ denotes the $K\times K$ diagonal matrix whose $(k,k)$-coefficient is $\lambda_k$. Thus, $\mathbf{\Sigma}$ and its eigenstructure are available in closed-form, and the latter coincides with the eigenstructure of $S^2_{\mathbf{x}}$ up to a scaling factor.

We let $\chi_{K}^2$ denote a chi-squared distribution with $K$ degrees of freedom, $\mathbb{P}(\chi_{K}^2\leq a)$ its cumulative distribution function (CDF) evaluated at $a\in\mathbb{R}$, and $q_{\chi_{K}^2}(p)$ its $p$-quantile for $p\in(0,1)$.

\section{Review of Rerandomization}
\label{sec: review}
{We follow the potential outcomes framework \citep{rubin1990comment,rubin2005causal}, where each unit $i$ has fixed potential outcomes $Y_i(1)$ and $Y_i(0)$, which denote the outcome for unit $i$ under treatment and control, respectively.} Thus, the observed outcome for unit $i$ is $y_i^{obs} = W_i Y_i(1) + (1 - W_i)Y_i(0)$. Define $\mathbf{y}^{obs} = (y_1^{obs} \,\dots\, y_N^{obs})^\top$ as the vector of observed outcomes. We focus on the average treatment effect as the causal estimand, defined as
\begin{align}
	\tau = \frac{1}{N} \sum_{i=1}^N [Y_i(1) - Y_i(0)] \label{eqn:ateEstimand}.
\end{align}
Furthermore, we focus on the mean-difference estimator
\begin{align}
	\hat{\tau} = \bar{\mathbf{y}}_T - \bar{\mathbf{y}}_C\,, \label{eqn:meanDiffEstimator}
\end{align}
where $\bar{\mathbf{y}}_T = N_T^{-1}\mathbf{W}^\top\mathbf{y}^{obs}$ and $\bar{\mathbf{y}}_C = N_C^{-1}(\mathbf{1}_N-\mathbf{W})^\top\mathbf{y}^{obs}$ are the average treatment and control outcomes, respectively. When conducting a randomized experiment, ideally we would like $\bar{\mathbf{x}}_T$ and $\bar{\mathbf{x}}_C$ to be close; otherwise, the estimator $\hat{\tau}$ could be confounded by imbalances in the covariate means.

\cite{morgan2012rerandomization} focused on a rerandomization scheme using the Mahalanobis distance to ensure that the covariate means are reasonably balanced for a particular treatment assignment. The Mahalanobis distance between the treatment and control covariate means is defined as
\begin{equation}
\begin{aligned}
	M &= (\bar{\mathbf{x}}_T - \bar{\mathbf{x}}_C)^\top \mathbf{\Sigma}^{-1}(\bar{\mathbf{x}}_T - \bar{\mathbf{x}}_C), \label{eqn:md}
\end{aligned}
\end{equation}
where the dependence of $M$ on the assignment vector $\mathbf{W}$ is implicit through $(\bar{\mathbf{x}}_T - \bar{\mathbf{x}}_C)$. \cite{morgan2012rerandomization} suggest randomizing units to treatment and control by performing independent draws from the distribution of $\mathbf{W}\,|\,\mathbf{x}$ until $M \leq a$ for some threshold $a \geq 0$. Hereafter, we refer to this procedure of randomizing units until $M \leq a$ as \textit{rerandomization}. The expected number draws until the first acceptable randomization is equal to $1/p_a$, where $p_a = \mathbb{P}(M \leq a\,|\,\mathbf{x})$ is the probability that a particular realization of $\mathbf{W}$ yields a Mahalanobis distance $M$ less than or equal to $a$. Thus, fixing $p_a$ effectively allocates an expected computational budget and induces a corresponding threshold $a$: the smaller the acceptance probability $p_a$, the smaller the threshold $a$ and thus the more balanced the two groups, but the larger the expected computational cost of drawing an acceptable $\mathbf{W}$. For example, to restrict rerandomization to the ``best'' 1\% randomizations, one would set $p_a = 0.01$, which implicitly sets $a$ equal to the $p_a$-quantile of the distribution of $M$ given $\mathbf{x}$. If one assumes $(\bar{\mathbf{x}}_T - \bar{\mathbf{x}}_C)\,|\,\mathbf{x} \sim \mathcal{N}(0, \mathbf{\Sigma})$, then $M\,|\,\mathbf{x} \sim \chi^2_K$, so that $a$ can be chosen equal to the $p_a$-quantile of a chi-squared distribution with $K$ degrees of freedom. The assumption $(\bar{\mathbf{x}}_T - \bar{\mathbf{x}}_C)\,|\,\mathbf{x}\sim\mathcal{N}(0,\mathbf{\Sigma})$ can be justified by invoking the finite population Central Limit Theorem \citep{erdos1959central,li2017general}. When the distribution of $M\,|\,\mathbf{x}$ is unknown, one can approximate it via Monte Carlo by simulating independent draws of $M\,|\,\mathbf{x}$ and setting $a$ to the $p_a$-quantile of $M$'s empirical distribution.

\cite{morgan2012rerandomization} established that the mean-difference estimator $\hat{\tau}$ under this rerandomization scheme is unbiased in estimating the average treatment effect $\tau$, i.e., that $\E{\hat{\tau}\,|\,\mathbf{x},M\leq a} = \tau$. Furthermore, they also established that under rerandomization, if $N_T=N_C$ and $(\bar{\mathbf{x}}_T - \bar{\mathbf{x}}_C)\,|\,\mathbf{x} \sim \mathcal{N}(0, \mathbf{\Sigma})$, then not only are the covariate mean differences centered at $0$, i.e., $\E{\bar{\mathbf{x}}_T - \bar{\mathbf{x}}_C\,|\,\mathbf{x},M\leq a} = 0$, but also they are more closely concentrated around $0$ than they would be under randomization. More precisely, \cite{morgan2012rerandomization} proved that
\setlength{\jot}{10pt}
\begin{align}
	\cov{\bar{\mathbf{x}}_T - \bar{\mathbf{x}}_C\,|\,\mathbf{x},M\leq a} &= v_a \cov{\bar{\mathbf{x}}_T - \bar{\mathbf{x}}_C\,|\,\mathbf{x}}, 
	\\
	\text{ with }\; v_a &= \frac{\mathbb{P}(\chi^2_{K+2} \leq a)}{\mathbb{P}(\chi^2_K \leq a)} \in(0,1). \label{eqn:va}
\end{align}
Therefore, under their assumptions, rerandomization using the Mahalanobis distance reduces the variance of each covariate mean difference by $100(1 - v_a)\%$ compared to randomization. \cite{morgan2012rerandomization} call this last property \textit{equally percent variance reducing} (EPVR). Thus, using the Mahalanobis distance for rerandomization can be quite appealing, but \cite{morgan2012rerandomization} rightly point out that non-EPVR rerandomization schemes may be preferable in settings with covariates of unequal importances. This is in part addressed by \cite{morgan2015rerandomization}, who developed a rerandomization scheme that incorporates tiers of covariates that vary in importance. However, this requires researchers to specify an explicit hierarchy of covariate importance, which may not be immediately clear, especially when the number of covariates is large. Furthermore, if there are strong collinearities amongst covariates such that $\mathbf{\Sigma}$ is degenerate and thus the $M$ in (\ref{eqn:md}) is undefined, then it is unclear how one should conduct the rerandomization scheme of \cite{morgan2012rerandomization} and its extensions \citep{morgan2015rerandomization,branson2016improving,li2018asymptotic,li2020rerandomization}.

\section{Ridge Rerandomization} \label{s:ridgeRerandomization}
As an alternative, we consider a modified Mahalanobis distance, defined as
\begin{equation}
M_{\lambda} = (\bar{\mathbf{x}}_T - \bar{\mathbf{x}}_C )^\top (\mathbf{\Sigma} + \lambda\; \mathbf{I}_K)^{-1} (\bar{\mathbf{x}}_T - \bar{\mathbf{x}}_C) \label{eqn:ridgeMD}
\end{equation}
for some prespecified $\lambda\geq 0$. Guidelines for choosing $\lambda$ will be provided in Section \ref{ss:guidelinesALambda}. The eigenvalues of $\mathbf{\Sigma}$ in \eqref{eqn:ridgeMD} are inflated in a way that is reminiscent of ridge regression \citep{hoerl1970ridge}. For this reason, we will refer to the quantity $M_{\lambda}$ as the \textit{ridge Mahalanobis distance}. To our knowledge, the ridge Mahalanobis distance has remained largely unexplored, except for \cite{kato1999handwritten}, who used it in an application for a Chinese and Japanese character recognition system. Our proposed rerandomization scheme, referred to as \textit{ridge rerandomization}, involves using the ridge Mahalanobis distance in place of the standard Mahalanobis distance within the rerandomization framework of \citet{morgan2012rerandomization}. In other words, one randomizes the assignment vector $\mathbf{W}$ until $M_{\lambda}\leq a_\lambda$ for some threshold $a_\lambda\geq0$.

In order to make a fair comparison between rerandomization and ridge rerandomization, we will fix the expected computational cost of ridge rerandomization by calibrating the respective thresholds so that 
\begin{equation}
	\mathbb{P}(M_\lambda \leq a_\lambda\,|\,\mathbf{x})=\mathbb{P}(M \leq a\,|\,\mathbf{x})=p_a.
	\label{eqn:computational_budget}
\end{equation}
Thus, fixing $p_a$ implicitly determines the pair $(\lambda,a_\lambda)$, so that for every fixed $\lambda\geq 0$ and $p_a\in(0,1)$ corresponds a unique $a_\lambda$ that satisfies \eqref{eqn:computational_budget}.

As we will discuss in Section \ref{ss:connections}, the ridge Mahalanobis distance alleviates collinearity among the covariate mean differences by placing higher importance on the directions that account for the most variation. In that section we also discuss how ridge rerandomization encapsulates a spectrum of other standard rerandomization schemes. But first, in Section \ref{ss:propertiesOfRidgeRerandomization} we establish several theoretical properties of ridge rerandomization for some prespecified $(\lambda, a_{\lambda})$, and in Section \ref{ss:guidelinesALambda} we provide guidelines for specifying $(\lambda, a_{\lambda})$. In Section \ref{ss:inference}, we discuss how to conduct inference for the average treatment effect $\tau$ after ridge rerandomization is used to design a randomized experiment.

\subsection{Properties of Ridge Rerandomization} \label{ss:propertiesOfRidgeRerandomization}

The following theorem establishes that, on average, the covariate means in the treatment and control groups are balanced under ridge rerandomization, and that $\hat{\tau}$ is an unbiased estimator of $\tau$ under ridge rerandomization.

\begin{theorem}[Unbiasedness under ridge rerandomization]
	\label{thm:covariateUnbiased} 
	Let $\lambda\geq 0$ and $a_{\lambda} \geq 0$ be some prespecified constants. If $N_T = N_C$, then
	\begin{equation*}
		\mathbb{E}[\bar{\mathbf{x}}_T - \bar{\mathbf{x}}_C \,|\,\mathbf{x}, M_{\lambda} \leq a_{\lambda}] = 0
	\end{equation*}
	and 
	\begin{equation*}
	\mathbb{E}[\hat{\tau} \,|\,\mathbf{x}, M_{\lambda} \leq a_{\lambda}] = \tau.
	\end{equation*}
\end{theorem}

\noindent
Theorem \ref{thm:covariateUnbiased} is a particular case of Theorem 2.1 and Corollary 2.2 from \cite{morgan2012rerandomization}. Theorem \ref{thm:covariateUnbiased} follows from the symmetry of $M_{\lambda}$ in treatment and control, in the sense that both assignments $\mathbf{W}$ and $(\mathbf{1}_N-\mathbf{W})$ yield the same value of $M_\lambda$. From \cite{morgan2012rerandomization}, we even have the stronger result that $\mathbb{E}[\bar{V}_T - \bar{V}_C \,|\,\mathbf{x}, M_{\lambda} \leq a_{\lambda}] = 0$ for any covariate $V$, regardless of whether $V$ is observed or not.

Now we establish the covariance structure of $(\bar{\mathbf{x}}_T - \bar{\mathbf{x}}_C)$ under ridge rerandomization. To do this, we first derive the exact distribution of $M_{\lambda}$. The following lemma establishes that if we assume $(\bar{\mathbf{x}}_T - \bar{\mathbf{x}}_C)\,|\,\mathbf{x}\sim\mathcal{N}(0,\mathbf{\Sigma})$, then $M_\lambda$ is distributed as a weighted sum of $K$ independent $\chi^2_1$ random variables, where the sizes of the weights are ordered in the same fashion as the sizes of the eigenvalues of $\mathbf{\Sigma}$.

\begin{lemma}[Distribution of $M_{\lambda}$]
\label{lemma:ridgeMDDistribution}
	Let $\lambda\geq 0$ be some prespecified constant. If $(\bar{\mathbf{x}}_T - \bar{\mathbf{x}}_C)\,|\,\mathbf{x}\sim\mathcal{N}(0,\mathbf{\Sigma})$, then
	\begin{equation}
\label{eq:weighted_sum}
M_{\lambda}\,|\,\mathbf{x} \;\;\sim\;\; \sum_{j=1}^K \frac{\lambda_j}{\lambda_j + \lambda} Z_j^2
\end{equation}
where $Z_1,...,Z_K \stackrel{\text{i.i.d.}}{\sim} N(0, 1)$ and $\lambda_1 \geq \dots \geq \lambda_K>0$ are the eigenvalues of $\mathbf{\Sigma}$.
\end{lemma}

\noindent
The proof of Lemma \ref{lemma:ridgeMDDistribution} is provided in the Appendix; see Section \ref{ss:proofLemmaRidgeMDDistribution}. Under the Normality assumption, the representation in \eqref{eq:weighted_sum} provides a straightforward way to simulate independent draws of $M_\lambda$, despite its CDF being typically intractable and requiring numerical approximations \citep[e.g., see][ and references therein]{bodenham2016comparison}.

Using Lemma \ref{lemma:ridgeMDDistribution}, we can derive the covariance structure of $\bar{\mathbf{x}}_T - \bar{\mathbf{x}}_C$ under ridge rerandomization, as stated by the following theorem.

\begin{theorem}[Covariance structure under ridge rerandomization]
	\label{thm:ridgeRerandomizationCovariance}
	Let $\lambda \geq 0$ and $a_{\lambda} \geq 0$ be some prespecified constants. If $(\bar{\mathbf{x}}_T - \bar{\mathbf{x}}_C)\,|\,\mathbf{x}\sim\mathcal{N}(0,\mathbf{\Sigma})$ and $N_T = N_C$, then
	\begin{equation}
	\emph{\text{Cov}}(\bar{\mathbf{x}}_T - \bar{\mathbf{x}}_C\,|\,\mathbf{x},M_{\lambda}\leq a_{\lambda}) = \boldsymbol{\Gamma} \emph{\textbf{Diag}}((\lambda_k\, d_{k,\lambda})_{1\leq k \leq K}) \boldsymbol{\Gamma}^\top
	\label{eqn:covariance_ridge}
	\end{equation}
	where $\boldsymbol{\Gamma}$ is the orthogonal matrix of eigenvectors of $\mathbf{\Sigma}$ corresponding to the ordered eigenvalues $\lambda_1 \geq \dots \geq \lambda_K > 0$, and for all $k=1,...,K$,
	\begin{equation}
	d_{k,\lambda} = \mathbb{E}\left[ Z_k^2 \,\left|\, \sum_{j=1}^K \frac{\lambda_j}{\lambda_j + \lambda} Z_j^2 \leq a_{\lambda} \right.\right]
	\label{eqn:diagonal_coeff_ridge}
	\end{equation}
	with $Z_1,...,Z_K \stackrel{\text{i.i.d.}}{\sim} N(0,1)$.
\end{theorem}
\noindent
The proof of Theorem \ref{thm:ridgeRerandomizationCovariance} is in the Appendix in Section \ref{ss:proofTheoremRidgeRerandomizationCovariance}. The quantities $d_{k,\lambda}$ are intractable functions of $\lambda$ and $a_{\lambda}$ and thus need to be approximated numerically, as explained in Section \ref{ss:calibration_alambda}. Conditioning on $M_\lambda\leq a_\lambda$ in \eqref{eqn:diagonal_coeff_ridge} effectively constrains the magnitude of the positive random variables $Z_k^2$. Since the weights $\lambda_k(\lambda_k+\lambda)^{-1}$ of their respective contributions to $M_\lambda$ are positive and non-increasing with $k=1,...,K$, we may conjecture that $0<d_{1,\lambda}\leq \dots \leq d_{K,\lambda}<1$. Possible directions for a proof may make use of Proposition 2.1 from \citet{palombi2013note} and Equation (A.1) from \citet{palombi2017numerical}.

Using the above results, we can now compare randomization, rerandomization, and ridge rerandomization. Under the assumptions stated in Theorem \ref{thm:ridgeRerandomizationCovariance}, the covariance matrices of $\bar{\mathbf{x}}_T - \bar{\mathbf{x}}_C$ under randomization, rerandomization, and ridge rerandomization can be respectively written as
\begin{align}
\cov{\bar{\mathbf{x}}_T - \bar{\mathbf{x}}_C\,|\,\mathbf{x}} &= \boldsymbol{\Gamma} \diag{(\lambda_k)_{1\leq k \leq K}} \boldsymbol{\Gamma}^\top, \label{eqn:randomizationVarianceReduction}
\\
\cov{\bar{\mathbf{x}}_T - \bar{\mathbf{x}}_C\,|\,\mathbf{x},M \leq a} &= \boldsymbol{\Gamma} \diag{(\lambda_k \,v_a)_{1\leq k \leq K}} \boldsymbol{\Gamma}^\top, \label{eqn:rerandomizationVarianceReduction}
\\
\cov{\bar{\mathbf{x}}_T - \bar{\mathbf{x}}_C\,|\,\mathbf{x},M_\lambda \leq a_\lambda} &= \boldsymbol{\Gamma} \diag{(\lambda_k \,d_{k,\lambda})_{1\leq k \leq K}} \boldsymbol{\Gamma}^\top. \label{eqn:ridgererandomizationVarianceReduction}
\end{align}
where \eqref{eqn:rerandomizationVarianceReduction} follows from Theorem 3.1 in \cite{morgan2012rerandomization} with $v_a\in(0,1)$, and \eqref{eqn:ridgererandomizationVarianceReduction} follows from Theorem \ref{thm:ridgeRerandomizationCovariance} with $d_{k,\lambda}\in(0,1)$ defined in \eqref{eqn:diagonal_coeff_ridge}. If we define new covariates $\mathbf{x}^*$ as the principal components of the original ones, i.e., $\mathbf{x}^* = \mathbf{x}\boldsymbol{\Gamma}$, then  \eqref{eqn:rerandomizationVarianceReduction} and \eqref{eqn:ridgererandomizationVarianceReduction} respectively yield
\begin{equation}
	\var{(\bar{\mathbf{x}}^*_{T} - \bar{\mathbf{x}}^*_{C})_k\,|\,\mathbf{x},M \leq a} = v_a \;\var{(\bar{\mathbf{x}}^*_{T} - \bar{\mathbf{x}}^*_{C})_k\,|\,\mathbf{x}}
	\label{eq:var_rerand_pc}
\end{equation}
and
\begin{equation}
\var{(\bar{\mathbf{x}}^*_{T} - \bar{\mathbf{x}}^*_{C})_k\,|\,\mathbf{x},M_\lambda \leq a_\lambda} = d_{k,\lambda} \;\var{(\bar{\mathbf{x}}^*_{T} - \bar{\mathbf{x}}^*_{C})_k\,|\,\mathbf{x}}
\label{eq:var_ridegrerand_pc}
\end{equation}
for all $k=1,...,K$, where $(\bar{\mathbf{x}}^*_{T} - \bar{\mathbf{x}}^*_{C})_k$ is the $k$-th \textit{principal component mean difference} between the treatment and control groups, i.e., the $k$-th coefficient of $\boldsymbol{\Gamma}^\top(\bar{\mathbf{x}}_{T}-\bar{\mathbf{x}}_{C})$. From \eqref{eq:var_rerand_pc} we see that rerandomization reduces the variances of the principal component mean differences equally by $100(1-v_a)\%$ and is thus EPVR for the principal components, as well as for the original covariates, as discussed in Section \ref{sec: review}. On the other hand, ridge rerandomization reduces these variances by unequal amounts: the variance of the $k$-th principal component mean difference is reduced by $100(1-d_{k,\lambda})\%$, {and because typically $0<d_{1,\lambda}\leq\dots\leq d_{K,\lambda}<1$, ridge rerandomization places more importance on the first principal components.}

Translating \eqref{eq:var_ridegrerand_pc} back to the original covariates yields the following corollary, which establishes that ridge rerandomization is always preferable over randomization in terms of reducing the variance of each covariate mean difference.
\begin{corollary}[Variance reduction for ridge rerandomization]
	\label{cor:varReductionRidgeRerandomization}
	Under the assumptions of Theorem \ref{thm:ridgeRerandomizationCovariance}, ridge rerandomization reduces the variance of the $k$-th covariate mean difference $(\bar{\mathbf{x}}_{T} - \bar{\mathbf{x}}_{C})_k$ by $100 \left(1 - v_{k,\lambda} \right)\%$, where
	\begin{align}
		v_{k,\lambda} = \frac{\left(\boldsymbol{\Gamma} \emph{\textbf{Diag}}\left((\lambda_j \,d_{j,\lambda})_{1\leq j \leq K}\right)\boldsymbol{\Gamma}^\top\right)_{kk}}{\mathbf{\Sigma}_{kk}} \label{eqn:vLambdaAK}
	\end{align}
	satisfies $v_{k,\lambda}\in(0,1)$, so that
	\begin{equation}
	\emph{\text{Var}}\left((\bar{\mathbf{x}}_{T} - \bar{\mathbf{x}}_{C})_k\,|\,\mathbf{x},M_\lambda \leq a_\lambda\right) \;<\; \emph{\text{Var}}\left((\bar{\mathbf{x}}_{T} - \bar{\mathbf{x}}_{C})_k\,|\,\mathbf{x}\right).
	\end{equation}
\end{corollary}

\noindent
The proof of Corollary \ref{cor:varReductionRidgeRerandomization} is provided in the Appendix; see Section \ref{ss:proofTheoremVarianceReductionBound}. Reducing the variance of the covariate mean differences is beneficial for precisely estimating the average treatment effect if the outcomes are correlated with the covariates. For example, Theorem 3.2 of \cite{morgan2012rerandomization} establishes that---under several assumptions, including additivity of the treatment effect---rerandomization reduces the variance of $\hat{\tau}$ defined in (\ref{eqn:meanDiffEstimator}) by $100(1 - v_a)R^2$ percent, where $R^2$ denotes the squared multiple correlation between the outcomes and the covariates. Now we establish how the variance of $\hat{\tau}$ behaves under ridge rerandomization.

In the rest of this section, we assume---as in \cite{morgan2012rerandomization}---that the treatment effect is additive. Without loss of generality, for all $i=1,...,N$, we can write the outcome of unit $i$ as
\begin{equation}
	Y_i(W_i) = \beta_0 + \mathbf{x}_{i\bullet}\boldsymbol{\beta} + \tau W_i + \epsilon_i
	\label{eq:regression_outcome}
\end{equation}
{where $\beta_0 + \mathbf{x}\boldsymbol{\beta}$ is the projection of the potential outcomes $\mathbf{Y}(0)=(Y_1(0) \dots Y_N(0))^\top$ onto the linear space spanned by $(\mathbf{1}, \mathbf{x})$, and $\epsilon_i\in\mathbb{R}$ captures any misspecification of the linear relationship between the outcomes and $\mathbf{x}$. Let $\bar{\boldsymbol{\epsilon}}_T = N_T^{-1}\mathbf{W}^\top\boldsymbol{\epsilon}$ and $\bar{\boldsymbol{\epsilon}}_C = N_C^{-1}(\mathbf{1}_N-\mathbf{W})^\top\boldsymbol{\epsilon}$, where $\boldsymbol{\epsilon} = (\epsilon_1 \dots \epsilon_N)^\top$.

Theorem \ref{thm:variancePercentReductionRidgeRerandomizationComparison} below establishes that the variance of $\hat{\tau}$ under ridge rerandomization is always less than or equal to the variance of $\hat{\tau}$ under randomization. Thus, ridge rerandomization always leads to a more precise treatment effect estimator than randomization.

\begin{theorem} \label{thm:variancePercentReductionRidgeRerandomizationComparison}
	Under the assumptions of Theorem \ref{thm:ridgeRerandomizationCovariance}, if $(\bar{\boldsymbol{\epsilon}}_T - \bar{\boldsymbol{\epsilon}}_C)$ is conditionally independent of $(\bar{\mathbf{x}}_T - \bar{\mathbf{x}}_C)$ given $\mathbf{x}$, and if there is an additive treatment effect, then
	\begin{equation}
		\emph{\text{Var}}(\hat{\tau}\,|\,\mathbf{x}) - \emph{\text{Var}}(\hat{\tau}\,|\,\mathbf{x},M_{\lambda} \leq a_{\lambda}) \;=\; \boldsymbol{\beta}^\top \boldsymbol{\Gamma} \emph{\textbf{Diag}} \left((\lambda_k \left(1 - d_{k,\lambda} \right))_{1\leq k \leq K} \right) \boldsymbol{\Gamma}^\top \boldsymbol{\beta}
		\nonumber
	\end{equation}
	so that we have
	\begin{align}
	\emph{\text{Var}}(\hat{\tau}\,|\,\mathbf{x},M_{\lambda} \leq a_{\lambda}) \;\leq\; \emph{\text{Var}}(\hat{\tau}\,|\,\mathbf{x})\,,
	\nonumber
	\end{align}
	where the equality holds if and only if $\boldsymbol{\beta}= \mathbf{0}_K$ in \eqref{eq:regression_outcome}.
\end{theorem}

\noindent
The proof of Theorem \ref{thm:variancePercentReductionRidgeRerandomizationComparison} is in the Appendix; see Section \ref{ss:proofTheoremVariancePercentReductionRidgeRerandomizationComparison}. The conditional independence assumption was also leveraged in the proof of Theorem 3.2 in \cite{morgan2012rerandomization}.

The fact that ridge rerandomization performs better than randomization is arguably a low bar, because this is the purpose of any rerandomization scheme. The following corollary quantifies how ridge rerandomization performs compared to the rerandomization scheme of \cite{morgan2012rerandomization}.

\begin{corollary}
\label{corr:varianceReductionRerandomizationRidgeRerandomizationComparison}
	Under the assumptions of Theorem \ref{thm:variancePercentReductionRidgeRerandomizationComparison}, the difference in variances of $\hat{\tau}$ between rerandomization and ridge rerandomization is
	\begin{align}
	\emph{\text{Var}}(\hat{\tau}\,|\,\mathbf{x},M \leq a) - \emph{\text{Var}}(\hat{\tau}\,|\,\mathbf{x},M_{\lambda} \leq a_{\lambda}) \;=\; \boldsymbol{\beta}^\top \boldsymbol{\Gamma} \emph{\textbf{Diag}} \left((\lambda_k \left(v_a - d_{k,\lambda} \right))_{1\leq k \leq K} \right) \boldsymbol{\Gamma}^\top \boldsymbol{\beta}.
	\nonumber
\end{align}
\end{corollary}
\noindent
It is not necessarily the case that $d_{k,\lambda} \leq v_a$ for all $k = 1,\dots,K$, and so it is not guaranteed that ridge rerandomization will perform better or worse than rerandomization in terms of treatment effect estimation. Ultimately, the comparison of rerandomization and ridge rerandomization depends on $\boldsymbol{\beta}$, which is typically not known until after the experiment has been conducted.

However, in Section \ref{ss:comparingTreatmentEffectEstimation}, we provide some heuristic arguments for when ridge rerandomization would be preferable over rerandomization, along with simulation evidence that confirms these heuristic arguments. In particular, we demonstrate that ridge rerandomization is preferable over rerandomization when there are strong collinearities among the covariates. We also discuss a ``worst-case scenario'' for ridge rerandomization, where $\boldsymbol{\beta}$ is specified such that ridge rerandomization should perform worse than rerandomization in terms of treatment effect estimation accuracy.

In order to implement ridge rerandomization, researchers must specify the threshold $a_{\lambda}\geq 0$ and the regularization parameter $\lambda\geq 0$. The next section provides guidelines for choosing these parameters.

\subsection{Guidelines for choosing $a_{\lambda}$ and $\lambda$} \label{ss:guidelinesALambda}
For ridge rerandomization, we recommend starting by specifying an acceptance probability $p_{a}\in(0,1)$, which then binds $\lambda$ and $a_\lambda$ together via the identity \eqref{eqn:computational_budget}. Once $p_a$ is fixed, there exists a uniquely determined threshold $a_\lambda\geq 0$ for each $\lambda\geq 0$ such that $\mathbb{P}(M_\lambda\leq a_\lambda\,|\,\mathbf{x}) = p_a$. As in \cite{morgan2012rerandomization}, acceptable treatment allocations under ridge rerandomization are generated by randomizing units to treatment and control until $M_{\lambda} \leq a_{\lambda}$. Thus, a smaller $p_{a_{\lambda}}$ leads to stronger covariate balance according to $M_{\lambda}$ at the expense of computation time.

The only choice that remains after fixing $p_a$ is the regularization parameter $\lambda\geq 0$. Section \ref{ss:calibration_alambda} details how $a_\lambda$ is automatically calibrated once we fix $p_a$ and $\lambda$. The choice of $\lambda$ is investigated in Section \ref{sss:choiceoflambda}, after discussing how to assess the performance of ridge rerandomization in Section \ref{ss:approximation_dk}.
\subsubsection{Calibration of $a_\lambda$} \label{ss:calibration_alambda}
Given $p_a$ and $\lambda$, we can choose to set $a_\lambda$ equal to the $p_a$-quantile of the quadratic form $Q_\lambda$ defined by
\begin{equation}
	Q_\lambda = \sum_{k=1}^K \frac{\lambda_k}{\lambda_k+\lambda} Z_k^2
	\label{eq:quadraticform}
\end{equation}
where $Z_1,...,Z_K\stackrel{\text{i.i.d.}}{\sim}\mathcal{N}(0,1)$. Such a choice of $a_\lambda$ is a good approximation of the $p_a$-quantile of $M_\lambda$, especially when $N$ is large enough for $(\bar{\mathbf{x}}_T-\bar{\mathbf{x}}_C)\,|\,\mathbf{x}$ to be approximately Normal, as motivated by Lemma \ref{lemma:ridgeMDDistribution}. Let $F_{Q_\lambda}(q) = \mathbb{P}(Q_\lambda\leq q)$ denote the CDF of $Q_\lambda$. Since $Q_\lambda$ is a weighted sum of independent $\chi_1^2$ variables, its characteristic function $\phi_{Q_\lambda}$ is given by $\phi_{Q_\lambda}(t) = \prod_{k=1}^{K}[1-2i\lambda_k(\lambda_k+\lambda)^{-1}t]^{-1/2}$, which can then be inverted to yield
\begin{equation*}
F_{Q_\lambda}(q) = \lim_{U\to+\infty}F_{Q_\lambda,U}(q)
\end{equation*}
where
\begin{equation}
	F_{Q_\lambda,U}(q)= \frac{1}{2} - \frac{1}{\pi}\int_{0}^{U}\;\frac{\sin\left(\frac{1}{2}\left[-t\,q + \sum_{k=1}^{K}\arctan\left(\frac{\lambda_k}{\lambda_k+\lambda}\,t\right)\right]\right)}{t\;\prod_{k=1}^{K}\left[1+\left(\frac{\lambda_k}{\lambda_k+\lambda}\right)^2t^2\right]^{1/4}}\;dt
\end{equation}
as detailed in Equation (3.2) of \citet{imhof1961computing}. In practice, for any fixed $U\geq 0$, $F_{Q_\lambda,U}(q)$ can be computed with arbitrary precision and at a negligible cost by using any (deterministic) univariate numerical integration scheme. We can then approximate $F_{Q_\lambda}(q)$ with $F_{Q_\lambda,U}(q)$ by choosing $U$ large enough. As explained in \citet{imhof1961computing}, the approximation tends to improve as the number of covariates $K$ increases, and one can guarantee a truncation error of at most $\xi>0$ in absolute value by choosing $U_\xi = [\xi\,\pi\,(K/2)\prod_{k=1}^{K}\sqrt{\lambda_k(\lambda_k+\lambda)^{-1}}]^{-2/K}$. Computationally cheaper but less accurate alternatives to approximate $F_{Q_\lambda}$ are discussed in \citet{bodenham2016comparison}.

Finally, we approximate the $p_a$-quantile of $Q_\lambda$ by
\begin{equation}
	\hat{a}_{\lambda} = \inf\{q\in\mathbb{R}:F_{Q_\lambda,U_\xi}(q)\geq p_a\}
	\label{eq:quantile_approx}
\end{equation}
i.e., the $p_a$-quantile of $F_{Q_\lambda,U}$. The hat on $\hat{a}_\lambda$ only reflects the distributional approximation of $M_\lambda$ by $Q_\lambda$, whereas the errors due to numerical integration and truncation can be regarded as virtually nonexistent compared to the Monte Carlo errors involved in the later approximations of $v_{k,\lambda}$. In the simulations of Section \ref{s:simulations}, we will use $\xi = 10^{-4}$ by default.

\subsubsection{Approximation of $d_{k,\lambda}$ and $v_{k,\lambda}$} \label{ss:approximation_dk}
We will use Corollary \ref{cor:varReductionRidgeRerandomization} and Theorem \ref{thm:ridgeRerandomizationCovariance} as a proxy for how ridge rerandomization improves the variance of each covariate mean difference as compared to rerandomization. We would like to set $(\lambda,a_\lambda)$ so that the $d_{k,\lambda}$'s defined in \eqref{eqn:diagonal_coeff_ridge} are small, in a sense to be made precise in the next section. To achieve this, we would need to compute $d_{k,\lambda}$ for all $k=1,...K$, which involves intractable conditional expectations. By considering $n$ simulated sets of $K$ independent variables $\widetilde{Z}_{ij}\stackrel{\text{i.i.d.}}{\sim}\mathcal{N}(0,1)$ for $i=1,..,n$ and $j=1,..,K$, the expectations appearing in \eqref{eqn:diagonal_coeff_ridge} can be consistently estimated via Monte Carlo, for all $k = 1,...,K$, by
\begin{equation}
	\hat{d}_{k,\lambda,n} = \frac{1}{\sum_{i=1}^n\mathds{1}(M_\lambda^{(i)}\,\leq\, \hat{a}_{\lambda})} \sum_{i=1}^n \widetilde{Z}_{ik}^{\,2}\,\mathds{1}(M_\lambda^{(i)}\,\leq\, \hat{a}_{\lambda})
	\label{eq:approxMonteCarlo_dk}
\end{equation}
with $M_\lambda^{(i)} = \sum_{j=1}^K \lambda_k(\lambda_k+\lambda)^{-1}\widetilde{Z}_{ij}^{\,2}$ and $\hat{a}_{\lambda}$ defined in \eqref{eq:quantile_approx}, where $\mathds{1}(A)$ denotes the indicator function of an event $A$. Using \eqref{eq:approxMonteCarlo_dk}, we can then estimate  $v_{k,\lambda}$ from Corollary \ref{cor:varReductionRidgeRerandomization} consistently as $n\to +\infty$, for all $k = 1,...,K$, by
\begin{equation}
\hat{v}_{k,\lambda,n} = \frac{\left(\boldsymbol{\Gamma} {\textbf{Diag}}\left((\lambda_j \,\hat{d}_{j,\lambda,n})_{1\leq j \leq K}\right)\boldsymbol{\Gamma}^\top\right)_{kk}}{\mathbf{\Sigma}_{kk}} \label{eqn:approx_vLambdaAK}
\end{equation}
For simplicity, we will regard the computational cost of generating $nK$ independent Normal variables as negligible compared to the expected cost of generating $1/p_a$ successive random assignment vectors and testing the acceptability of each assignment, since the former can be done in parallel at virtually the same cost as generating one single Normal random variable.

\subsubsection{Choosing $\lambda$} \label{sss:choiceoflambda}
In this section, assume that $p_a$ has been fixed. Note that choosing $\lambda = 0$ corresponds to rerandomization using the Mahalanobis distance. Thus, we would only choose some $\lambda > 0$ if it is preferable over rerandomization, in the following sense. There are many metrics that could be used for comparing rerandomization and ridge rerandomization; for simplicity, we focus on the average percent reduction in variance across covariate mean differences. Arguably, one rerandomization scheme is preferable over another if it is able to achieve a higher average reduction in variance across covariates. Thus, ideally, we would only choose a particular $\lambda>0$ if $K^{-1}\sum_{k=1}^K v_{k,\lambda}< v_a$. In practice, we will use the criterion
\begin{align}
	\frac{1}{K} \sum_{k=1}^K \hat{v}_{k,\lambda,n} \;<\; v_a   \label{eqn:ridgeLambdaInequality}
\end{align}
\noindent
where $v_a$ and $\hat{v}_{k,\lambda,n}$ are respectively defined in \eqref{eqn:va} and \eqref{eqn:approx_vLambdaAK}, with $a$ being set to $q_{\chi_{K}^2}(p_a)$, i.e., the choice of $a$ as recommended by \cite{morgan2012rerandomization}. Proving the existence of some $\lambda>0$ such that \eqref{eqn:ridgeLambdaInequality} holds is challenging, so we propose the following iterative procedure for choosing such a $\lambda>0$ if it exists. The procedure relies on \eqref{eqn:va}, \eqref{eq:quantile_approx}, and \eqref{eqn:approx_vLambdaAK}, where the auxiliary Normal variables $\widetilde{Z}_{ij}$ only need to be simulated once and can then be reused when testing different values of $\lambda$.
\\\\
\noindent\fbox{%
\parbox{\textwidth}{%
\vspace*{0.2cm}
\textbf{ Procedure for finding a desirable $\lambda\geq 0$}
\begin{enumerate}
	\item Specify $p_a\in(0,1)$, $n\geq 1$, $\delta>0$, and $\varepsilon>0$.
	\item Initialize $\lambda = 0$ and $\Lambda = \emptyset$.
	\item While $|(\lambda+\delta)\hat{a}_{\lambda+\delta} - \lambda\hat{a}_{\lambda}|>\varepsilon$:
	\begin{itemize}
		\item Set $\lambda = \lambda + \delta$.
		\item If $\;\displaystyle\frac{1}{K} \sum_{k=1}^K \hat{v}_{k,\lambda,n} \;<\; \frac{\mathbb{P}(\chi_{K+2}^2\;\leq\; q_{\chi_{K}^2}(p_a))}{p_a}$, then set $\Lambda =  \Lambda\cup\{\lambda\}$.
	\end{itemize}
	\item If $\Lambda=\emptyset$, then return $\lambda = 0$. \\
	Else, define $c_k = \lambda_k^2\,(\sum_{j=1}^{K}\lambda_j^2)^{-1}$ for all $k=1,...,K$, and return:
	\begin{equation}
		\lambda_\star = \operatorname*{argmin}_{\lambda\in\Lambda}\left(\;\sum_{k=1}^{K}c_k\,\hat{d}^{\,2}_{k,\lambda,n}-\left(\sum_{k=1}^{K}c_k\,\hat{d}_{k,\lambda,n}\right)^2\;\right).
		\label{eq:best_lambda}
	\end{equation}
\end{enumerate}
}}
\\\\\\
The justification of our proposed procedure stems from the following facts. By definition, we have $\mathbb{P}(M_\lambda\leq a_\lambda\,|\,\mathbf{x})=p_a$ for all $\lambda\geq 0$. By taking the limit as $\lambda\to+\infty$ under the assumptions of Lemma \ref{lemma:ridgeMDDistribution},  we get
\begin{equation*}
	p_a = \lim_{\lambda\to+\infty}\mathbb{P}\left(\sum_{k=1}^{K}\frac{\lambda_k}{\lambda_k+\lambda}Z_k^2\leq a_\lambda\right) = \lim_{\lambda\to+\infty}\mathbb{P}\left(\sum_{k=1}^{K}\lambda_kZ_k^2\leq \lambda \,a_\lambda\right)
\end{equation*}
so that 
\begin{equation}
\lambda\,a_\lambda\quad\xrightarrow[\lambda\to+\infty]{}\quad q^*(p_a)
\label{eq:limit_lambda_a}
\end{equation}
where $q^*(p_a)$ is the $p_a$-quantile of the distribution of $\sum_{k=1}^{K}\lambda_kZ_k^2$. This in turn implies that, for all $k=1,...,K$, we have
\begin{equation}
	v_{k,\lambda}\quad\xrightarrow[\lambda\to+\infty]{}\quad \frac{\left(\boldsymbol{\Gamma} {\textbf{Diag}}\left((\lambda_j \,d^*_{j})_{1\leq j \leq K}\right)\boldsymbol{\Gamma}^\top\right)_{kk}}{\mathbf{\Sigma}_{kk}}
	\label{eq:limit_vk}
\end{equation}
where $d_k^* = \E{Z_k^2|\sum_{k=1}^{K}\lambda_kZ_k^2\leq q^*(p_a)}$ for all $k=1,...,K$. Since the limits in \eqref{eq:limit_vk} are strictly positive, this shows that increasing $\lambda$ beyond a certain value will no longer yield any practical gain. This is in line with the intuition that the ridge Mahalanobis distance degenerates to the Euclidean distance when $\lambda\to+\infty$, as discussed further in Section \ref{ss:connections}. Thus, in practice, it is sufficient to search for $\lambda$ only over a bounded range of values. The lower bound $\lambda=0$ corresponds to rerandomization with the standard Mahalanobis distance; the upper bound is determined dynamically via Step 3, which is guaranteed to stop in finite time by using an argument similar to \eqref{eq:limit_lambda_a}. The step size $\delta$ can be chosen as a fraction of the smallest strictly positive gap between consecutive eigenvalues, i.e., $\min\{\lambda_k-\lambda_{k-1}:k=1,...,K \text{ such that }\lambda_k>\lambda_{k-1}\}$ with the convention $\lambda_0=0$. Finally, among all the acceptable $\lambda$'s satisfying \eqref{eqn:ridgeLambdaInequality}, Step 4 returns the $\lambda_\star$ that aims at altering the conditional covariance structure of $(\bar{\mathbf{x}}_T-\bar{\mathbf{x}}_C)$ the least, in the sense of minimizing the distance between $\cov{\bar{\mathbf{x}}_T-\bar{\mathbf{x}}_C|\mathbf{x},M_\lambda\leq \hat{a}_\lambda}$ and the linear span of $\mathbf{\Sigma}$, i.e.,
\begin{equation*}
	\lambda_\star = \operatorname*{argmin}_{\lambda\in\Lambda}\left(\min_{c\in\mathbb{R}}\,{\left\|c\mathbf{\Sigma}-\boldsymbol{\Gamma} {\textbf{Diag}}\left((\lambda_j \,\hat{d}_{j,\lambda,n})_{1\leq j \leq K}\right)\boldsymbol{\Gamma}^\top\right\|}\;\right)
\end{equation*}
where $\|\mathbf{\Sigma}\| = \sqrt{\text{tr}(\mathbf{\Sigma}^\top \mathbf{\Sigma})} = \sum_{k=1}^K \lambda_k^2$ stands for the Frobenius norm. The inner minimization can be written as
\begin{equation*}
	\min_{c\in\mathbb{R}}\left(\sum_{k=1}^{K}\lambda_k^2\left(c-\hat{d}_{k,\lambda,n}\right)^2\right)
\end{equation*}
which is attained at $c_\star = \sum_{k=1}^{K}c_k\,\hat{d}_{k,\lambda,n}$ with $c_k = \lambda_k^2\,(\sum_{j=1}^{K}\lambda_j^2)^{-1}$ for all $k=1,...,K$, thus yielding \eqref{eq:best_lambda}. The outer minimization is then straightforward since the set $\Lambda$ of candidates is finite by construction. We elaborate on why we choose a $\lambda_\star$ that alters the conditional covariance structure of $(\bar{\mathbf{x}}_T-\bar{\mathbf{x}}_C)$ the least in Section \ref{ss:connections}.

When the set $\Lambda$ is empty, we simply return $\lambda=0$ (which corresponds to typical rerandomization), although the following heuristic argument illustrates why we would expect the existence of at least one $\lambda$ such that (\ref{eqn:ridgeLambdaInequality}) holds. The rerandomization scheme of \cite{morgan2012rerandomization} spreads the benefits of variance reduction across all $K$ covariates equally; however, note that the term $v_a = {\mathbb{P}(\chi^2_{K+2} \leq q_{\chi_{K}^2}(p_a))}/p_a$ is monotonically increasing in the number of covariates $K$ for a fixed acceptance probability $p_a$. A consequence of this is that if one can instead determine a smaller set of $K_e<K$ covariates that is most relevant, then that smaller set of covariates can benefit from a greater variance reduction than what would be achieved by considering all $K$ covariates. As we mentioned at the end of Section \ref{sec: review}, this idea was partially addressed in \cite{morgan2015rerandomization}, which extended the rerandomization scheme of \cite{morgan2012rerandomization} to allow for tiers of covariate importance specified by the researcher, such that the most important covariates receive the most variance reduction. Ridge rerandomization, on the other hand, automatically specifies a hierarchy of importance based on the eigenstructure of the covariate mean differences. To provide intuition for this idea, consider a simple case where the smallest $(K - K_e)$ eigenvalues $\lambda_{K_e+1},...,\lambda_K$ are all arbitrarily close to $0$. In this case, we can find $\lambda>0$ such that ${\lambda_j}({\lambda_j+\lambda})^{-1}\approx 1$ for the $K_e$ largest eigenvalues and ${\lambda_j}({\lambda_j+\lambda})^{-1}\approx 0$ for the remaining $K - K_e$ eigenvalues, so that $M_{\lambda}$ would be approximately distributed as $\chi_{k_{e}}^2$ with an effective number of degrees of freedom $K_e$ strictly less than $K$. For some fixed acceptance probability $p_a\in(0,1)$ and corresponding thresholds $a_e = q_{\chi_{K_e}^2}(p_a)$ and $a = q_{\chi_{K}^2}(p_a)$, we would then have
\begin{align}
 v_{a_e} = \frac{\mathbb{P}(\chi^2_{K_e+2} \leq q_{\chi_{K_e}^2}(p_a))}{p_a} \;<\; \frac{\mathbb{P}(\chi^2_{K+2} \leq q_{\chi_{K}^2}(p_a))}{p_a} = v_a
 \end{align}
 since $p_a$ is fixed and $K_e < K$. The relative variance reduction for ridge rerandomization would then be $(1 - v_{a_e})$ for the first $K_e$ principal components---which in this simple example make up the total variation in the covariate mean differences---while the relative variance reduction for rerandomization would be $(1 - v_a)<(1-v_{a_e})$ for the $K$ covariates. Thus, in this case, ridge rerandomization would achieve a greater variance reduction on a lower-dimensional representation of the covariates than typical rerandomization.

 This heuristic argument also hints that our method has connections to a principal-components rerandomization scheme, where one instead balances on some lower dimension of principal components rather than on the covariates themselves. We discuss this point further in Section \ref{ss:connections}.

\subsection{Connections to Other Rerandomization Schemes} \label{ss:connections}

Ridge rerandomization has connections to other rerandomization schemes. Ridge rerandomization requires specifying the parameter $\lambda$; thus, consider two extreme choices of $\lambda$:
\begin{enumerate}
	\item $\lambda = 0$: $M_{\lambda} = M$, i.e., $M_{\lambda}$ corresponds to the Mahalanobis distance.
	\item $\lambda \to+ \infty$: $M_\lambda \approx \lambda^{-1}\|\bar{\mathbf{x}}_T-\bar{\mathbf{x}}_C\|^2$, i.e., $M_{\lambda}$ tends to a scaled Euclidean distance.
\end{enumerate}
Thus, for any finite $\lambda > 0$, the distance defined by $M_{\lambda}$ can be regarded as a compromise between the Mahalanobis and Euclidean distances. Rerandomization using the Euclidean distance is similar to a rerandomization scheme that places a separate caliper on each covariate, which was proposed by \cite{moulton2004covariate}, \cite{maclure2006measuring}, \cite{bruhn2009pursuit}, and \cite{cox2009randomization}. However, \cite{morgan2012rerandomization} note that such a rerandomization scheme is not affinely invariant and does not preserve the correlation structure of $(\bar{\mathbf{x}}_T - \bar{\mathbf{x}}_C)$ across acceptable randomizations. See \cite{morgan2012rerandomization} for a full discussion of the benefits of using affinely invariant rerandomization criteria. As discussed in Section \ref{sss:choiceoflambda}, our proposed procedure aims for larger variance reductions of important covariate mean differences while mitigating the perturbation of the correlation structure of $(\bar{\mathbf{x}}_T - \bar{\mathbf{x}}_C)$.

As an illustration, consider a randomized experiment where $N_T = N_C = 50$ units are assigned to treatment and control; and furthermore, where there are two correlated covariates, generated as $x_{1j} \stackrel{\text{i.i.d.}}{\sim} N(0,1)$ and $x_{2j} \stackrel{\text{i.i.d.}}{\sim} N(x_{1i}, 1)$ for $j = 1,...,N$. Figure \ref{fig:covariateMeanDifferenceDistributionPlot} shows the distribution of $(\bar{\mathbf{x}}_T - \bar{\mathbf{x}}_C)\,|\,\mathbf{x}$ across 1000 randomizations, rerandomizations (with $p_a = 0.1$), ridge rerandomizations (with $p_a = 0.1$ and $\lambda = 0.005$), and rerandomizations using the Euclidean distance instead of the Mahalanobis distance.

All three rerandomization schemes reduce the variance of $(\bar{\mathbf{x}}_T - \bar{\mathbf{x}}_C)_k\,|\,\mathbf{x}$ for $k\in\{1,2\}$, compared to randomization; however, rerandomization using the Euclidean distance destroys the correlation structure of $(\bar{\mathbf{x}}_T - \bar{\mathbf{x}}_C)\,|\,\mathbf{x}$, while rerandomization and ridge rerandomization largely maintain it. This provides further motivation for Step 4 of the procedure presented in Section \ref{sss:choiceoflambda}.

\begin{figure}[H]
	\centering
	\includegraphics[width=\linewidth]{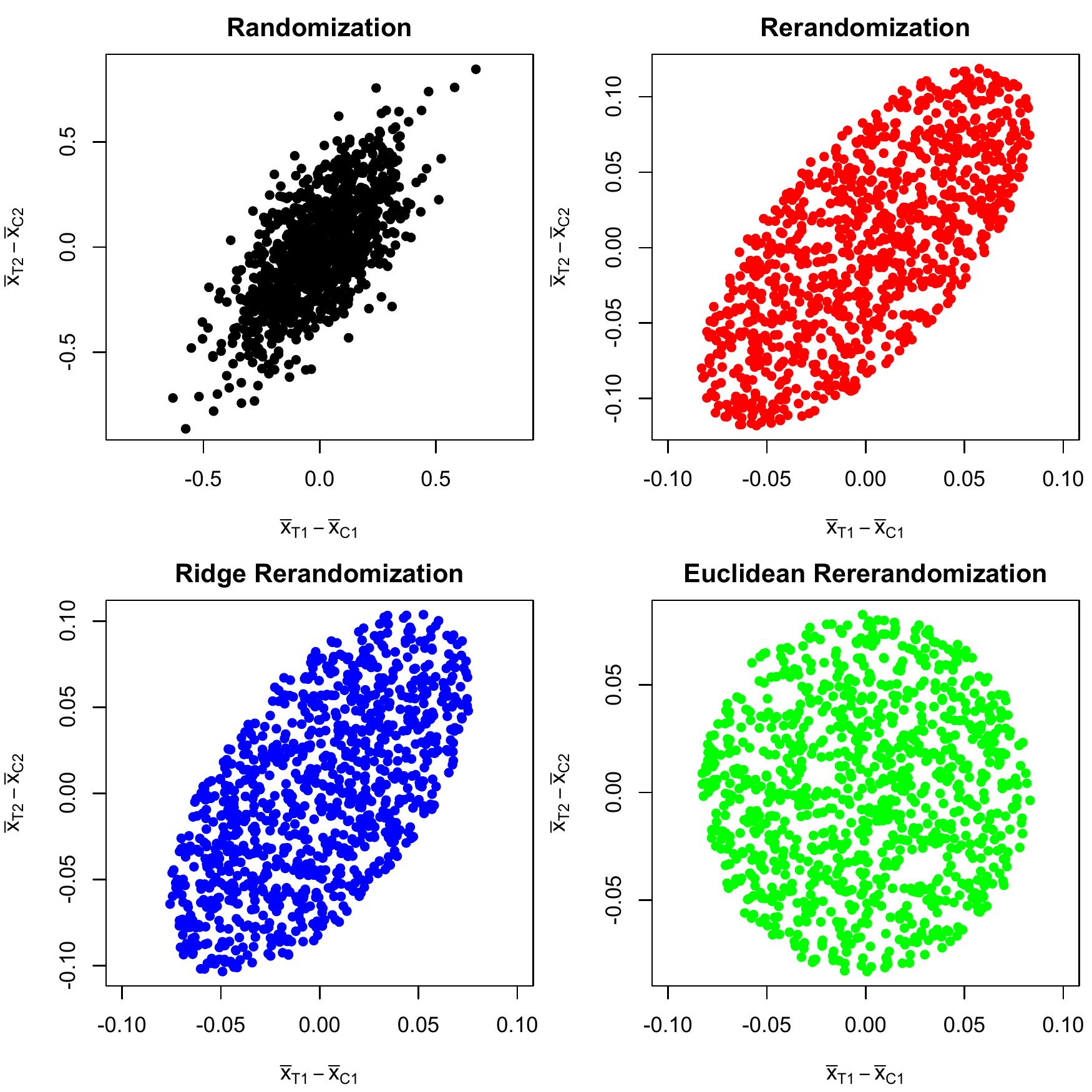}
	\caption{Distribution of $(\bar{\mathbf{x}}_T - \bar{\mathbf{x}}_C)\,|\,\mathbf{x}$ under randomization, rerandomization (with $p_a = 0.1$), ridge rerandomization (with $p_a = 0.1$ and $\lambda = 0.005$), and rerandomization using the Euclidean distance. Note the difference in scale for the randomization plot for ease of comparison.}
	\label{fig:covariateMeanDifferenceDistributionPlot}
\end{figure}

Furthermore, as discussed in Sections \ref{ss:propertiesOfRidgeRerandomization} and \ref{sss:choiceoflambda}, ridge rerandomization can be regarded as a ``soft-thresholding" version of a rerandomization scheme that would focus solely on the first $K_e<K$ principal components of $(\bar{\mathbf{x}}_T - \bar{\mathbf{x}}_C)$. A ``hard-thresholding" rerandomization scheme would use a truncated version $M_{K_e}$ of the Mahalanobis distance, defined as
\begin{equation*}
	M_{K_e} = (\bar{\mathbf{x}}_T - \bar{\mathbf{x}}_C)^\top \hat{\boldsymbol{\Sigma}}_{K_e}^{-1} (\bar{\mathbf{x}}_T - \bar{\mathbf{x}}_C)
\end{equation*}
with 
\begin{equation*}
	\boldsymbol{\Sigma}_{K_e} = \boldsymbol{\Gamma}\diag{(\lambda_1,...,\lambda_{K_e},0,...,0)}\boldsymbol{\Gamma}^\top
\end{equation*}
i.e., $\boldsymbol{\Sigma}_{K_e}$ artificially sets the smallest $(K-K_e)$ eigenvalues of $\boldsymbol{\Sigma}$ to $0$. This scheme would then be EPVR for the first $K_e$ principal components of $(\bar{\mathbf{x}}_T - \bar{\mathbf{x}}_C)$---although not necessarily EPVR for the original covariates themselves---but would effectively ignore the components associated with the smallest $(K-K_e)$ eigenvalues of $\boldsymbol{\Sigma}$.

Therefore, ridge rerandomization is a flexible experimental design strategy that encapsulates a class of rerandomization schemes, thus making it worth further investigation in future work. We expand on this point in Section \ref{s:discussion}.

\subsection{Conducting Inference After Ridge Rerandomization} \label{ss:inference}

Here we outline how to conduct inference for the average treatment effect after ridge rerandomization has been used to conduct an experiment. In general, there are Neymanian, Bayesian, and randomization-based modes of inference for analyzing randomized experiments \citep{imbens2015causal}. The Neymanian mode of inference relies on asymptotic approximations for the variance of the mean-difference estimator $\hat{\tau}$; such results are well-established for completely randomized experiments \citep{splawa1990application}, paired experiments \citep{imai2008variance}, blocked experiments \citep{miratrix2013adjusting,pashley2017insights}, and randomized experiments with stages of random sampling \citep{branson2019sampling}. In a seminal paper, \cite{li2018asymptotic} derived many asymptotic results for rerandomized experiments (as discussed in \cite{morgan2012rerandomization}), thereby establishing Neymanian inference for such experiments. The results therein rely on various properties of the Mahalanobis distance, which---as established by our results---differ from the properties of the ridge Mahalanobis distance. As a consequence, the theory developed in \cite{li2018asymptotic} cannot be readily applied to ridge rerandomized experiments, and a promising line of future work is deriving asymptotic results for ridge rerandomized experiments. Asymptotic results could also be used to establish Bayesian inference for such experiments, which would be particularly useful given that one's preference for rerandomization or ridge rerandomization may depend on their prior knowledge of $\boldsymbol{\beta}$, as suggested by Corollary \ref{corr:varianceReductionRerandomizationRidgeRerandomizationComparison}. Addressing these complications is beyond the scope of this paper. Instead, we focus on randomization-based inference, because it can be readily applied to ridge rerandomization.

Randomization-based inference focuses on inverting sharp null hypotheses that define the relationship between the potential outcomes in terms of treatment effects. The most common null hypothesis is that of an additive treatment effect $\tau$, such that the hypothesis $H_0^{\tau}: Y_i(1) = Y_i(0) + \tau$ holds for all $i=1,\dots,N$. Confidence intervals derived from inverting this hypothesis were first established by \cite{hodges1963estimates} and have since been popularized for analyzing randomized experiments (e.g., see \cite{rosenbaum2002overt} and \cite{imbens2015causal}). Here we briefly review how to obtain randomization-based confidence intervals for completely randomized experiments, and then we extend them to ridge rerandomized experiments.

As first proposed by \cite{hodges1963estimates}, a valid randomization-based confidence interval is the set of $\tau$ such that we fail to reject $H_0^{\tau}$; such inversion of a hypothesis is a classical way to obtain a confidence set \citep{kempthorne1969behaviour}. To obtain a valid $p$-value for $H_0^{\tau}$, a key insight is that, if $H_0^{\tau}$ holds, then one has full knowledge of the potential outcomes for all units: If we observe the outcome under control for a particular unit, we know that the outcome under treatment for that unit is simply the observed outcome plus $\tau$. As a result, for any hypothetical randomization, a test statistic---such as the mean difference estimator, $\hat{\tau}$---can be computed. To obtain a $p$-value for $H_0^{\tau}$ under randomization, one follows this simple three-step procedure:
\begin{enumerate}
	\item Generate many hypothetical randomizations, $\mathbf{w}^{(1)},\dots,\mathbf{w}^{(M)}$, by permuting the observed treatment indicator.
	\item Compute a test statistic $t(\mathbf{w}, \mathbf{x}, \mathbf{y})$, such as the mean-difference estimator, across the randomizations $\mathbf{w}^{(1)},\dots,\mathbf{w}^{(M)}$ assuming $H_0^{\tau}$ is true.
	\item Compute the randomization-based $p$-value, defined as
	\begin{align}
		p = \frac{1 + \sum_{m=1}^M \mathds{1} \left( |t(\mathbf{w}^{(m)}, \mathbf{x}, \mathbf{y})| > |t^{obs}| \right)}{M+1}
	\end{align}
\end{enumerate}
where $t^{obs}$ is the observed test statistic and $\mathds{1}(\cdot)$ denotes the indicator function. The additional 1 in the numerator and the denominator induces a very small amount of bias in order to validly control the Type 1 error rate and is a standard correction for randomization test $p$-values \citep{phipson2010permutation}. Modern statistical software allows one to readily invert $H_0^{\tau}$ after Step 1 is completed (in Section \ref{s:simulations}, we will use the \texttt{R} package \texttt{ri} \citep{aronow2012ri} to do this), thereby producing randomization-based confidence intervals. This makes the extension to ridge rerandomization quite straightforward: In Step 1, one generates many hypothetical \textit{ridge} rerandomizations (instead of randomizations), and then proceeds as usual to conduct randomization-based inference. This is identical to the approach discussed in \cite{morgan2012rerandomization} for obtaining confidence intervals under rerandomization, except using hypothetical ridge rerandomizations instead of hypothetical rerandomizations. This can also be viewed as inverting a \textit{conditional} randomization test, where we condition on the fact that the ridge rerandomization balance criterion has been fulfilled \citep{hennessy2016conditional,branson2019randomization}. As we shall see in Section \ref{s:simulations}, confidence intervals for ridge rerandomized experiments are much more precise than intervals for completely randomized experiments, and often more precise than intervals for rerandomized experiments, especially in high dimensional and/or collinearity settings.

\section{Simulations} \label{s:simulations}

We now provide simulation evidence that supports the hueristic argument presented in Section \ref{ss:guidelinesALambda} and suggests when ridge rerandomization is an effective experimental design strategy. First, we will consider conducting an experiment where covariates are linearly related with the outcome, and then we will consider when they are nonlinearly related. Throughout, we will compare rerandomization and ridge rerandomization in terms of (1) their ability to balance covariates, (2) their ability to produce precise treatment effect estimators, and (3) their ability to produce precise confidence intervals. We find that ridge rerandomization is particularly preferable over rerandomization in high-dimensional or high-collinearity settings.

\subsection{Simulation Setup} \label{ss:simulationSetup}

Consider $N = 100$ units, 50 of which are to be assigned to treatment and 50 are to be assigned to control. Let $\mathbf{x}$ be a $N \times K$ covariate matrix, generated as
\begin{align}
	\mathbf{x} \sim \mathcal{N} \left( \begin{pmatrix}
		0 \\
		\vdots \\
		0
	\end{pmatrix},
	\begin{pmatrix}
		1 & \rho & \cdots & \rho \\
		\rho & 1 & \cdots & \rho \\
		\vdots & \vdots & \ddots & \vdots \\
		\rho & \rho & \cdots & 1
	\end{pmatrix} \label{eqn:manyCovariatesModel}
	\right)
\end{align}
where $0 \leq \rho < 1$. The parameter $\rho$ corresponds to the correlation among the covariates. Furthermore, let $Y_i(1)$ and $Y_i(0)$ be the potential outcomes under treatment and control, respectively, for unit $i$, generated as
\begin{equation}
\begin{aligned}
	Y_i(0) &\sim N \left( \mathbf{x} \boldsymbol{\beta}, 1 \right) \label{eqn:manyCovariatesPotentialOutcomesModel} \\
	Y_i(1) &= Y_i(0) + \tau
\end{aligned}
\end{equation}
For this simulation study, we set the treatment effect to be $\tau = 1$. In the above, the covariates are linearly related with the outcome; we conduct additional simulations where covariates are nonlinearly related with the outcome in Section \ref{ss:simsNonlinear}. Across simulations, we consider number of covariates $K \in \{10,\dots,90\}$ and correlation parameter $\rho \in \{0, 0.1,\dots,0.9\}$. We discuss choices for $\boldsymbol{\beta}$ in Section \ref{ss:comparingTreatmentEffectEstimation}. We also considered data-generating processes where covariances varied among covariates and where there are an uneven number of units assigned to treatment and control (i.e., unbalanced designs). However, the results for these other scenarios were largely the same as those for the above data-generating process, and so for ease of exposition we focus on results for the case where the covariates are generated from (\ref{eqn:manyCovariatesModel}) and the potential outcomes are generated from (\ref{eqn:manyCovariatesPotentialOutcomesModel}).

We will consider three experimental design strategies for assigning units to treatment and control:
\begin{enumerate}
	\item \textbf{Randomization}: Randomize $50$ units to treatment and $50$ to control.
	\item \textbf{Rerandomization}: Randomize $50$ units to treatment and $50$ to control until $M \leq a$, where $M$ is the Mahalanobis distance defined in (\ref{eqn:md}).
	\item \textbf{Ridge Rerandomization}: Randomize $50$ units to treatment and $50$ to control until $M_{\lambda} \leq a_{\lambda}$, where $M_{\lambda}$ is the ridge Mahalanobis distance defined in (\ref{eqn:ridgeMD}).
\end{enumerate}

For each choice of $K$, $\rho$, and $\boldsymbol{\beta}$, we ran randomization, rerandomization, and ridge rerandomization 1000 times. For rerandomization and ridge rerandomization, we set $p_a = 0.1$, which corresponds to randomizing within the 10\% ``best'' randomizations according to the Mahalanobis distance and ridge Mahalanobis distance, respectively. Furthermore, for ridge rerandomization, we used the procedure in Section \ref{sss:choiceoflambda} for selecting $\lambda$, with $n = 1000$, $\delta = 0.01$, and $\epsilon = 10^{-4}$. The value $\lambda = 0.01$ was selected for most $K$ and $\rho$, and occasionally $\lambda = 0.02$ was selected.

First, in Section \ref{ss:comparingCovariateBalance}, we compare how these three methods balanced the covariates $\mathbf{x}$, and so the $\boldsymbol{\beta}$ parameter in (\ref{eqn:manyCovariatesPotentialOutcomesModel}) is irrelevant for this section. Then, in Section \ref{ss:comparingTreatmentEffectEstimation}, we compare the accuracy of treatment effect estimators and precision of confidence intervals for each method; in this case, the specification of $\boldsymbol{\beta}$ is consequential.

\subsection{Comparing Covariate Balance Across Randomizations} \label{ss:comparingCovariateBalance}

First, we computed the covariate mean differences across each randomization, rerandomization, and ridge rerandomization. Figure \ref{fig:ridgeVSrerandPlotVarReduction} shows how much rerandomization and ridge rerandomization reduced the variance of $\bar{\mathbf{x}}_T - \bar{\mathbf{x}}_C$ (averaged across covariates) compared to randomization for data generated from (\ref{eqn:manyCovariatesModel}). For rerandomization, the average variance reduction decreases as $K$ increases (an observation previously made in \cite{morgan2012rerandomization}), and it stays largely the same across values of $\rho$ for fixed $K$. As for ridge rerandomization, the average variance reduction also decreases as $K$ increases, but the average variance reduction increases as $\rho$ increases, i.e., as there is more collinearity in $\mathbf{x}$. Finally, the right-hand plot in Figure \ref{fig:ridgeVSrerandPlotVarReduction} shows that ridge rerandomization has a higher average variance reduction than rerandomization; furthermore, the advantage of ridge rerandomization over rerandomization increases in both $K$ and $\rho$. This suggests that ridge rerandomization may be particularly preferable over rerandomization in the presence of many covariates and/or high collinearity among covariates, which is intuitive given the motivation of ridge regression \citep{hoerl1970ridge}.

\begin{figure}[H]
	\centerline{\includegraphics[width=\linewidth]{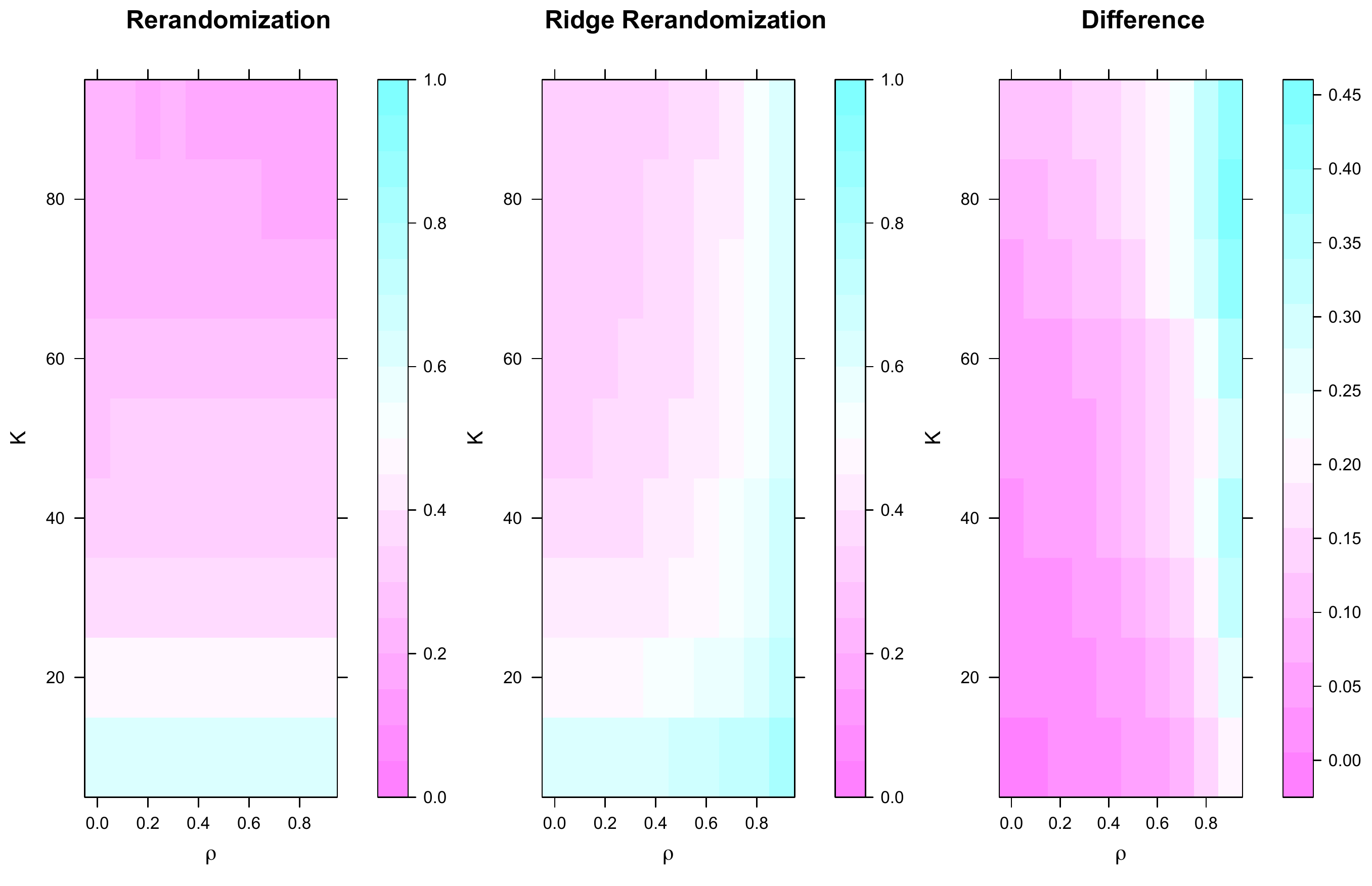}}
	\caption{Variance reduction averaged across covariates for rerandomization and ridge rerandomization, as well as their difference (ridge rerandomization minus rerandomization, i.e., the second plot minus the first).}
	\label{fig:ridgeVSrerandPlotVarReduction}
\end{figure}

\subsection{Comparing Treatment Effect Estimation Accuracy Across Randomizations} \label{ss:comparingTreatmentEffectEstimation}

Reducing the variance of each covariate mean difference leads to more precise treatment effect estimates if the covariates are related to the outcome, as in (\ref{eqn:manyCovariatesPotentialOutcomesModel}). The extent to which the covariates are related to the outcome depends on the $\boldsymbol{\beta}$ parameter. Theorem \ref{thm:variancePercentReductionRidgeRerandomizationComparison} guarantees that ridge rerandomization will improve inference for the average treatment effect, compared to randomization, regardless of $\boldsymbol{\beta}$. However, Corollary \ref{corr:varianceReductionRerandomizationRidgeRerandomizationComparison} establishes that $\boldsymbol{\beta}$ dictates whether rerandomization or ridge rerandomization will perform better in terms of treatment effect estimation accuracy. First we will consider a $\boldsymbol{\beta}$ where the covariates are equally related to the outcome, and in this case ridge rerandomization performs better than rerandomization. Then, we will consider a $\boldsymbol{\beta}$ which---according to our theoretical results---should put ridge rerandomization in the worst light as compared to rerandomization.

\subsubsection{One Choice of $\boldsymbol{\beta}$}

Consider $\boldsymbol{\beta} = \mathbf{1}_K$. Because the covariates have been standardized to have the same scale, such a $\boldsymbol{\beta}$ implies that all of the covariates are equally important in affecting the outcome. For each of the 1000 randomizations, rerandomizations, and ridge rerandomizations generated for each $K \in \{10,\dots,90\}$ and $\rho \in \{0, 0.1,\dots,0.9\}$, we computed the mean-difference estimator $\hat{\tau}$. Then, we computed the MSE of $\hat{\tau}$ across the 1000 randomizations, rerandomizations, and ridge rerandomizations for each $K$ and $\rho$. Figure \ref{fig:ridgeVSrerandTreatmentEffectMSE} shows the MSE of rerandomization and ridge rerandomization relative to the MSE of randomization. A lower relative MSE represents a more accurate treatment effect estimator, compared to how that estimator would behave under randomization.

\begin{figure}[H]
	\centerline{\includegraphics[width=\linewidth]{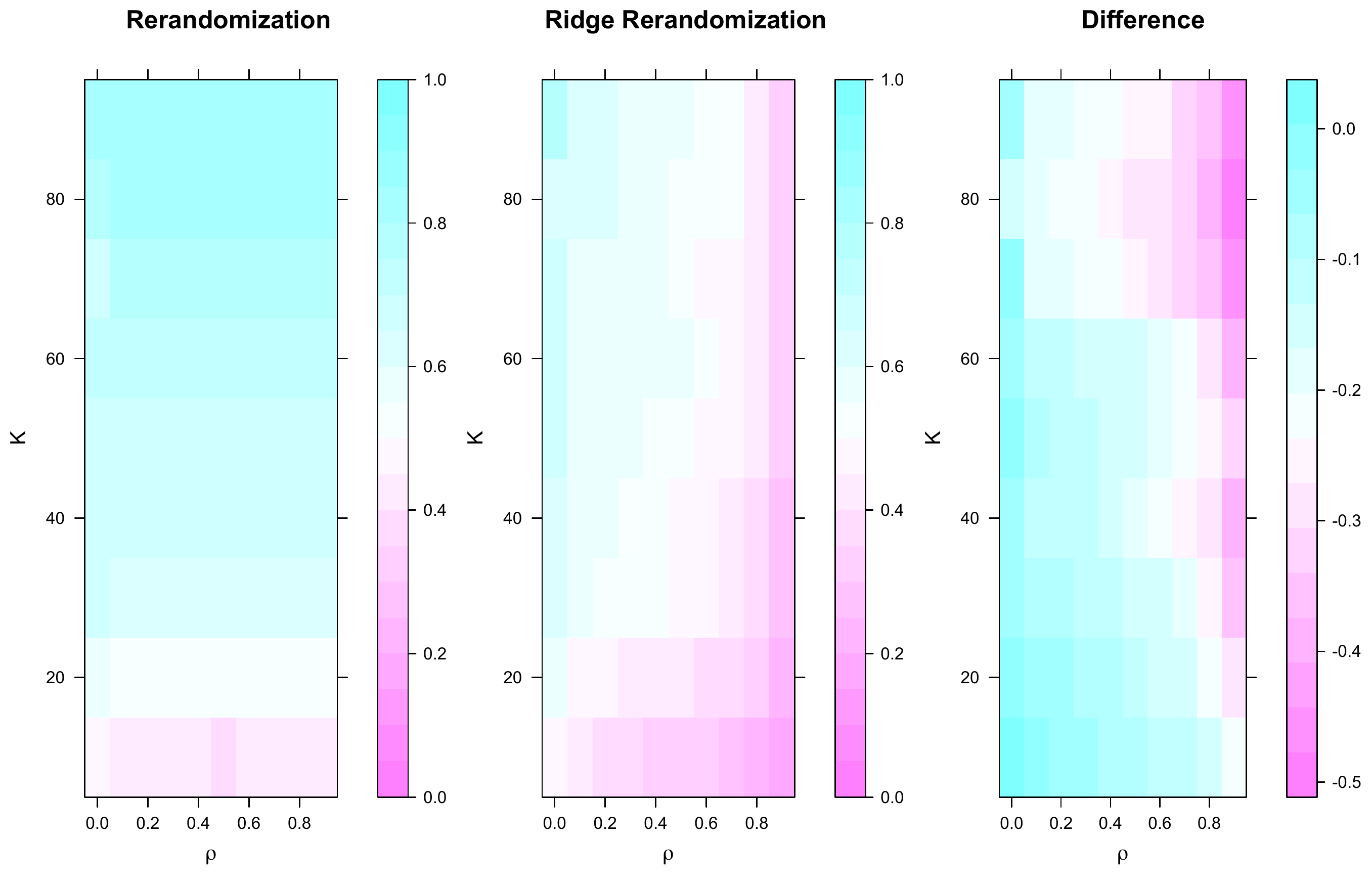}}
	\caption{Relative MSE of $\hat{\tau} = \bar{y}_T - \bar{y}_C$ under rerandomization and ridge rerandomization (relative to randomization) when $\boldsymbol{\beta} = \mathbf{1}_K$ in (\ref{eqn:manyCovariatesPotentialOutcomesModel}), as well as the difference in relative MSE between the two (i.e., the second plot minus the first).}
	\label{fig:ridgeVSrerandTreatmentEffectMSE}
\end{figure}

Three observations can be made about Figure \ref{fig:ridgeVSrerandTreatmentEffectMSE}. First, both rerandomization and ridge rerandomization reduce the MSE of $\hat{\tau}$ compared to randomization: the relative MSE for both methods is always less than 1. Second, for rerandomization, the relative MSE stays constant across values of $\rho$ and decreases as $K$ decreases. Meanwhile, for ridge rerandomization, the relative MSE decreases as $\rho$ increases and $K$ decreases. Third, for this choice of $\boldsymbol{\beta}$, ridge rerandomization reduces the MSE of the treatment effect estimator more so than rerandomization, especially when $K$ and/or $\rho$ is large. These last two observations reflect the variance reduction behavior observed in Figure \ref{fig:ridgeVSrerandPlotVarReduction}.

Meanwhile, for each randomization, rerandomization, and ridge rerandomization, we generated a 95\% confidence interval for the average treatment effect using the procedure outlined in Section \ref{ss:inference}. Regardless of the procedure used, coverage was near 95\%. This is unsurprising, because these intervals were constructed by inverting randomization tests that are valid for their corresponding assignment mechanism; see \cite{edgington2007randomization} and \cite{good2013permutation} for classical results on the validity of randomization tests. However, the width of these intervals differed across these three procedures: Figure \ref{fig:ridgeVSrerandCI} compares the relative average interval length (compared to randomization) for rerandomization and ridge rerandomization. For the first two plots in Figure \ref{fig:ridgeVSrerandCI}, a number closer to 1 indicates intervals that are closer in length to intervals under randomization. Meanwhile, for the right-most plot in Figure \ref{fig:ridgeVSrerandCI}, a more negative number indicates more narrow confidence intervals for ridge rerandomization, as compared to rerandomization. The qualitative results are identical to the previous results: Ridge rerandomization tends to provide narrower confidence intervals as the covariates' dimension and/or collinearity increases.

\begin{figure}[H]
	\centerline{\includegraphics[width=\linewidth]{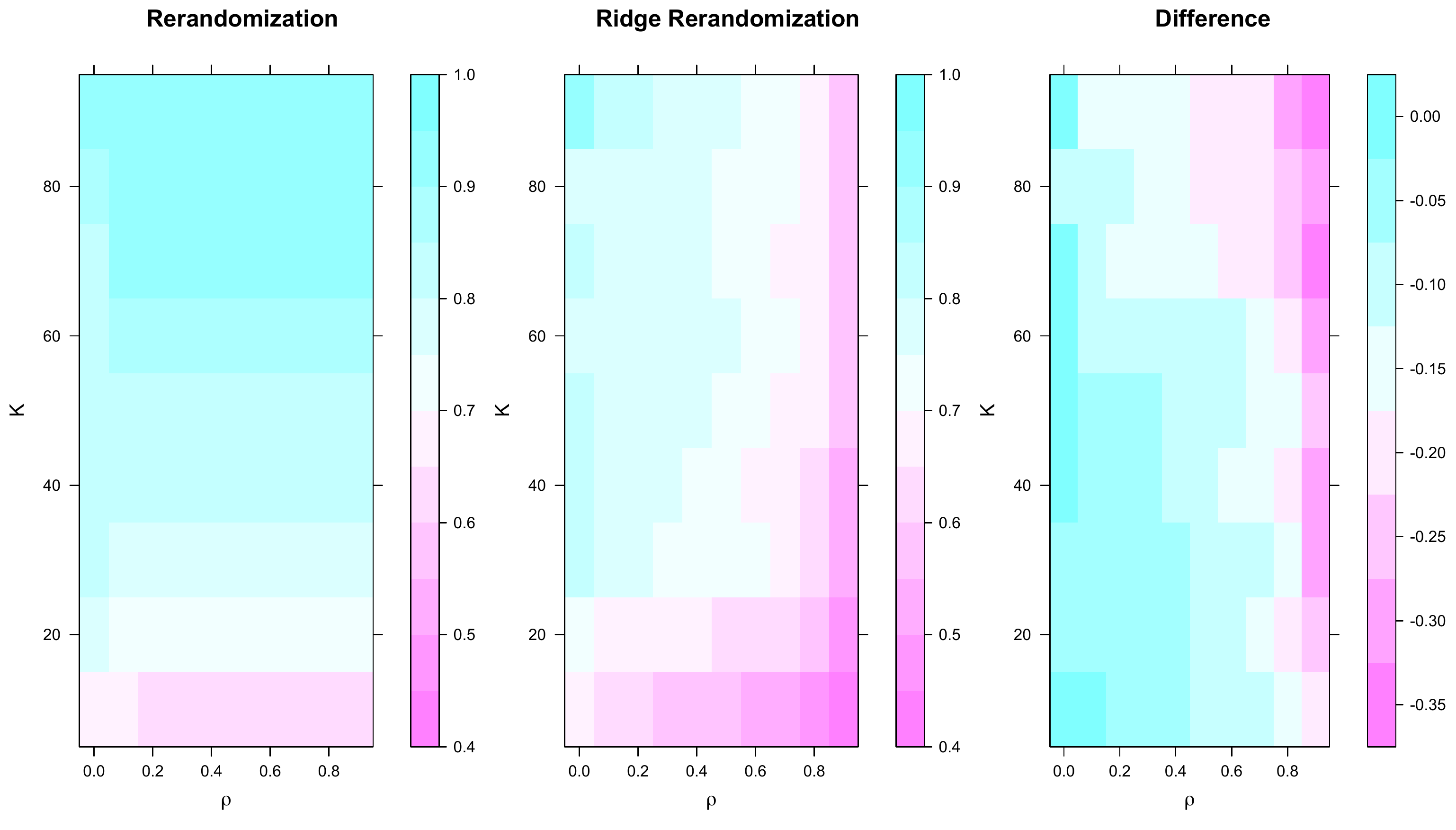}}
	\caption{Relative average confidence interval width under rerandomization and ridge rerandomization (relative to randomization) when $\boldsymbol{\beta} = \mathbf{1}_K$ in (\ref{eqn:manyCovariatesPotentialOutcomesModel}), as well as the difference between the two (i.e., the second plot minus the first).}
	\label{fig:ridgeVSrerandCI}
\end{figure}

\subsubsection{A Choice of $\boldsymbol{\beta}$ where Ridge Rerandomization has the Least Competitive Advantage over Rerandomization} \label{sss:badBeta}

As can be seen by Corollary \ref{corr:varianceReductionRerandomizationRidgeRerandomizationComparison}, there may exist $\boldsymbol{\beta}$ where rerandomization performs better than ridge rerandomization. {To assess how poorly ridge rerandomization can perform compared to rerandomization,} now we will specify a $\boldsymbol{\beta}$ that puts ridge rerandomization in the worst light when comparing it to rerandomization in terms of treatment effect estimation accuracy.

{Under the assumptions of Corollary \ref{corr:varianceReductionRerandomizationRidgeRerandomizationComparison}, the difference in treatment effect estimation accuracy between rerandomization and ridge rerandomization is given by $\Delta = \boldsymbol{\beta}^\top \boldsymbol{\Gamma} {\textbf{Diag}} \left((\lambda_k \left(v_a - d_{k,\lambda} \right))_{1\leq k \leq K} \right) \boldsymbol{\Gamma}^\top \boldsymbol{\beta}$, which can be artificially minimized with respect to $\boldsymbol{\beta}$, subject to some constraint on $\boldsymbol{\beta}$ for the minimum to exist, e.g., $\|\boldsymbol{\beta}\|\leq 1$. If $d_{k,\lambda}<v_a$ for all $k=1,...,K$, then ridge rerandomization dominates rerandomization since $\Delta>0$ for all $\boldsymbol{\beta}\neq 0$, and these schemes are only tied when $\Delta=0$ for $\boldsymbol{\beta}=0$, i.e., the covariates are uncorrelated with the outcomes. In other cases, we can define $\boldsymbol{\beta}^* = \boldsymbol{\Gamma}_{\bullet k^*}$ where $\boldsymbol{\Gamma}_{\bullet k^*}$ is the $k^*$-th column of $\boldsymbol{\Gamma}$ and $k^*=\operatorname{argmin}_{1\leq k \leq K}(v_a - d_{k,\lambda})$. We would typically have $k^*=K$ when the $d_{k,\lambda}$'s are non-increasing. By construction, $\boldsymbol{\beta}^*$ minimizes $\Delta$ over $\{\boldsymbol{\beta}\in\mathbb{R}^K:\|\boldsymbol{\beta}\|\leq 1\}$ and yields $\Delta<0$ as negative as possible. This is equivalent to $\boldsymbol{\beta}$ being in the direction that accounts for the least variation in the covariates. While such a case is unlikely, we consider such a $\boldsymbol{\beta}$ to see how much worse ridge rerandomization performs as compared to rerandomization in this scenario.}

Figure \ref{fig:ridgeVSrerandTreatmentEffectMSEBadBeta} shows the relative MSE (as compared to randomization) for rerandomization and ridge rerandomization for this specification of $\boldsymbol{\beta}$. Interestingly, there are occasions where rerandomization and ridge rerandomization have relative MSEs greater than 1, i.e., when they perform worse than randomization in terms of treatment effect estimation accuracy. At first this may be surprising, especially when findings from \cite{morgan2012rerandomization} guarantee that rerandomization should perform better than randomization. However, in this case, $\boldsymbol{\beta}$ is in the direction of the last principal component of the covariate space, meaning that the covariates have nearly no relationship with the outcomes. Thus, the relative MSE that we see in the first two plots of Figure \ref{fig:ridgeVSrerandTreatmentEffectMSEBadBeta} is more or less the behavior we would expect if we compared 1000 randomizations to 1000 other randomizations. Furthermore, from the third plot in Figure \ref{fig:ridgeVSrerandTreatmentEffectMSEBadBeta}, we can see that rerandomization occasionally performs better than ridge rerandomization---particularly when $K$ is small---but the differences in relative MSE across simulations are somewhat centered around zero. Meanwhile, Figure \ref{fig:ridgeVSrerandCIBadBeta} compares the relative average confidence interval length for rerandomization and ridge rerandomization, and the qualitative results are largely the same as the relative MSE results: Rerandomization and ridge rerandomization are fairly comparable, but rerandomization tends to provide slightly narrower confidence intervals for low-dimensional covariates.

Note that this specification of $\boldsymbol{\beta}$ is a unit vector. We could have scaled $\boldsymbol{\beta}$ arbitrarily large, and, as a result, the differences in the last plots of Figure \ref{fig:ridgeVSrerandTreatmentEffectMSEBadBeta} and \ref{fig:ridgeVSrerandCIBadBeta} could have been made arbitrarily large. Thus, ridge rerandomization can perform much worse than rerandomization when $\boldsymbol{\beta}$ exhibits particularly large effects in the direction of the last principal component of the covariate space, especially when the number of covariates is small. Practically speaking, such a scenario is unlikely, but it is a scenario that researchers should acknowledge and consider when comparing rerandomization and ridge rerandomization.

\begin{figure}[H]
	\centerline{\includegraphics[width=\linewidth]{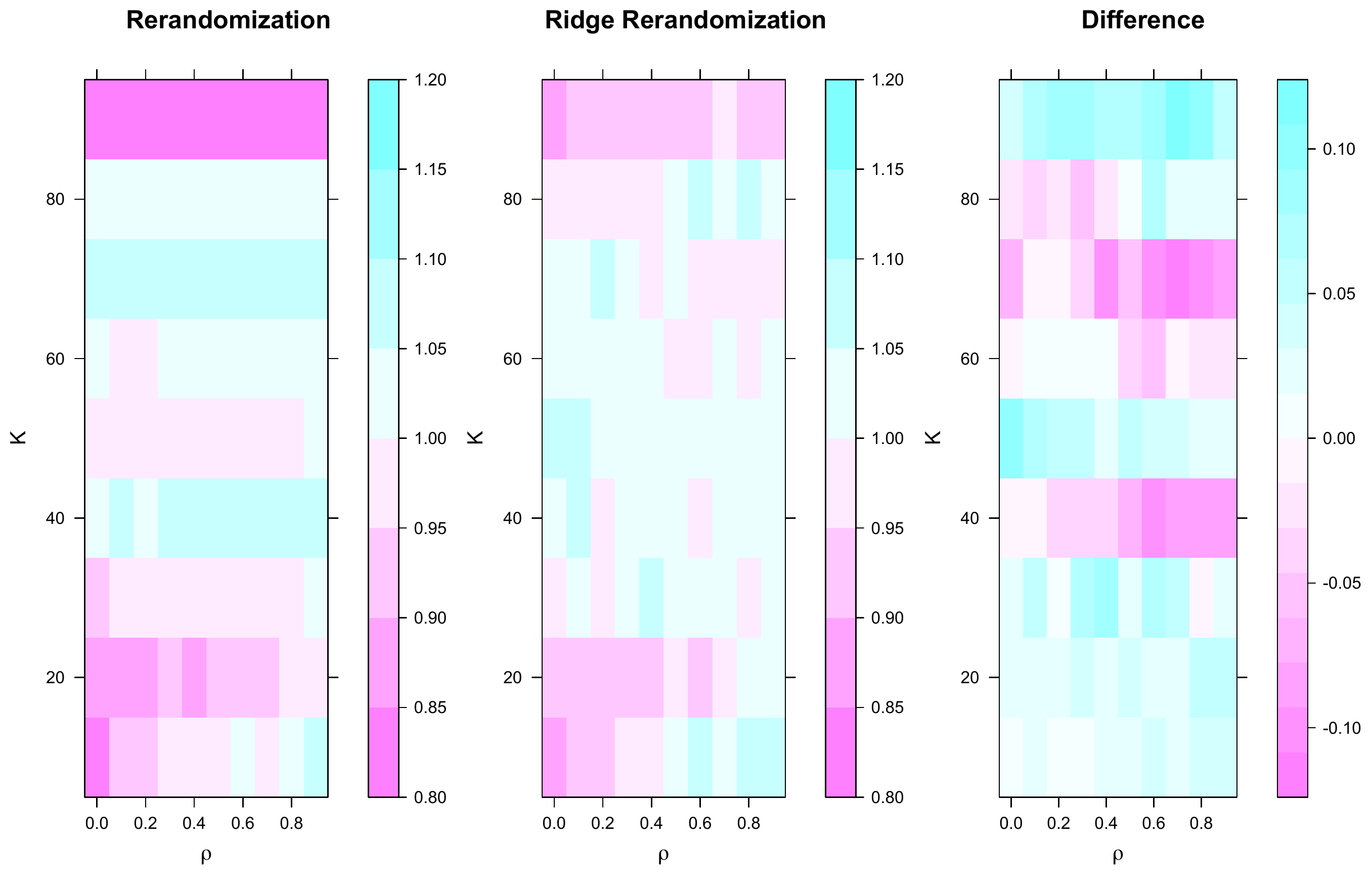}}
	\caption{Relative MSE of $\hat{\tau} = \bar{y}_T - \bar{y}_C$ under rerandomization and ridge rerandomization (relative to randomization) for the $\boldsymbol{\beta}$ such that ridge rerandomization has the least competitive advantage over rerandomization, as well as the difference in relative MSE between the two (i.e., the second plot minus the first).}
	\label{fig:ridgeVSrerandTreatmentEffectMSEBadBeta}
\end{figure}

\begin{figure}[H]
	\centerline{\includegraphics[width=\linewidth]{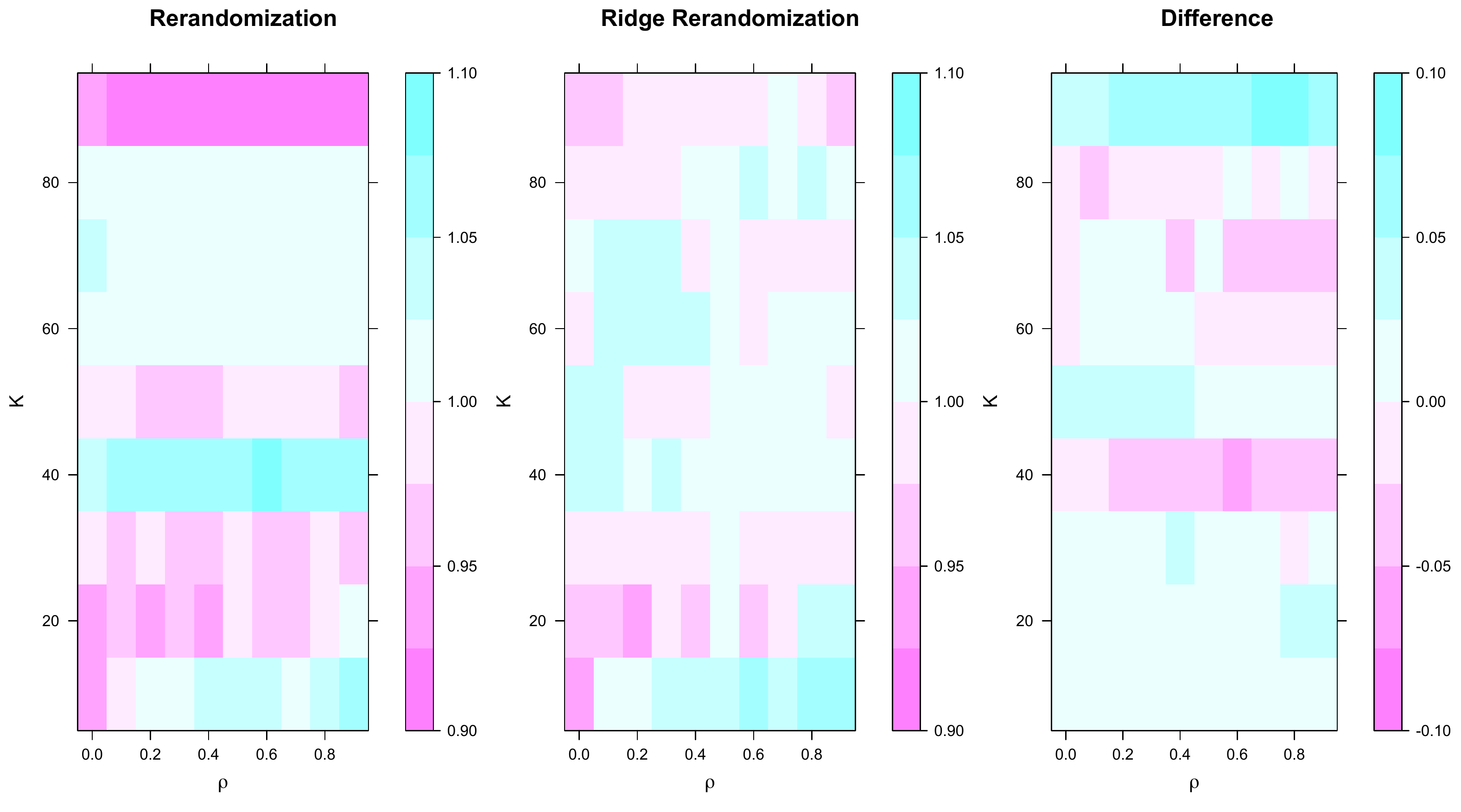}}
	\caption{Relative average confidence interval width under rerandomization and ridge rerandomization (relative to randomization) for the $\boldsymbol{\beta}$ such that ridge rerandomization has the least competitive advantage over rerandomization, as well as the difference between the two (i.e., the second plot minus the first).}
	\label{fig:ridgeVSrerandCIBadBeta}
\end{figure}

\subsection{Additional Simulations: Considering Nonlinearity} \label{ss:simsNonlinear}

In the above simulations, we generated the potential outcomes to be linearly related with the covariates. This may be considered a ``well-specified'' case, because ridge rerandomization only aims to balance the first moments of the covariates. In this section, we consider the same simulation setup as Section \ref{ss:simulationSetup}, except instead of using (\ref{eqn:manyCovariatesPotentialOutcomesModel}) to generate the potential outcomes, we use the following model: 
\begin{equation}
\begin{aligned}
	Y_i(0) &\sim N \left( \exp(\mathbf{x}) \boldsymbol{\beta}, 1 \right) \label{eqn:manyCovariatesPotentialOutcomesMisspecifiedModel} \\
	Y_i(1) &= Y_i(0) + \tau
\end{aligned}
\end{equation}
where $\exp(\mathbf{x})$ denotes the matrix of values $e^{\mathbf{x}}$. Again we set $\tau = 1$ and $\boldsymbol{\beta} = \mathbf{1}_K$ and consider $K \in \{10,\dots,90\}$ and $\rho \in \{0, 0.1,\dots,0.9\}$. This alternative model does not affect rerandomization and ridge rerandomization's ability to balance covariates, but it does affect their ability to precisely estimate treatment effects. Figure \ref{fig:ridgeVSrerandTreatmentEffectMSEMisspecified} compares the relative MSE (compared to randomization) of rerandomization and ridge rerandomization, and Figure \ref{fig:ridgeVSrerandCIMisspecified} does the same for relative average confidence interval length. Although ridge rerandomization does not have as clear of an advantage over rerandomization in this misspecified scenario, it still tends to perform better than rerandomization in high-dimensional and high-collinearity settings. Furthermore, both rerandomization and ridge rerandomization still provide more precise inference for the average treatment effect compared to randomization, although not as much as when the potential outcomes were generated from a linear model. This is because the covariates still have some linear relationship with the covariates, and thus one can still obtain more precise estimators and intervals for the average treatment effect by balancing the first moments of the covariates \citep{li2018asymptotic}. In short, the results presented here are largely the same as the previous results, where the potential outcomes were linearly related with the covariates. 

\begin{figure}[H]
	\centerline{\includegraphics[width=\linewidth]{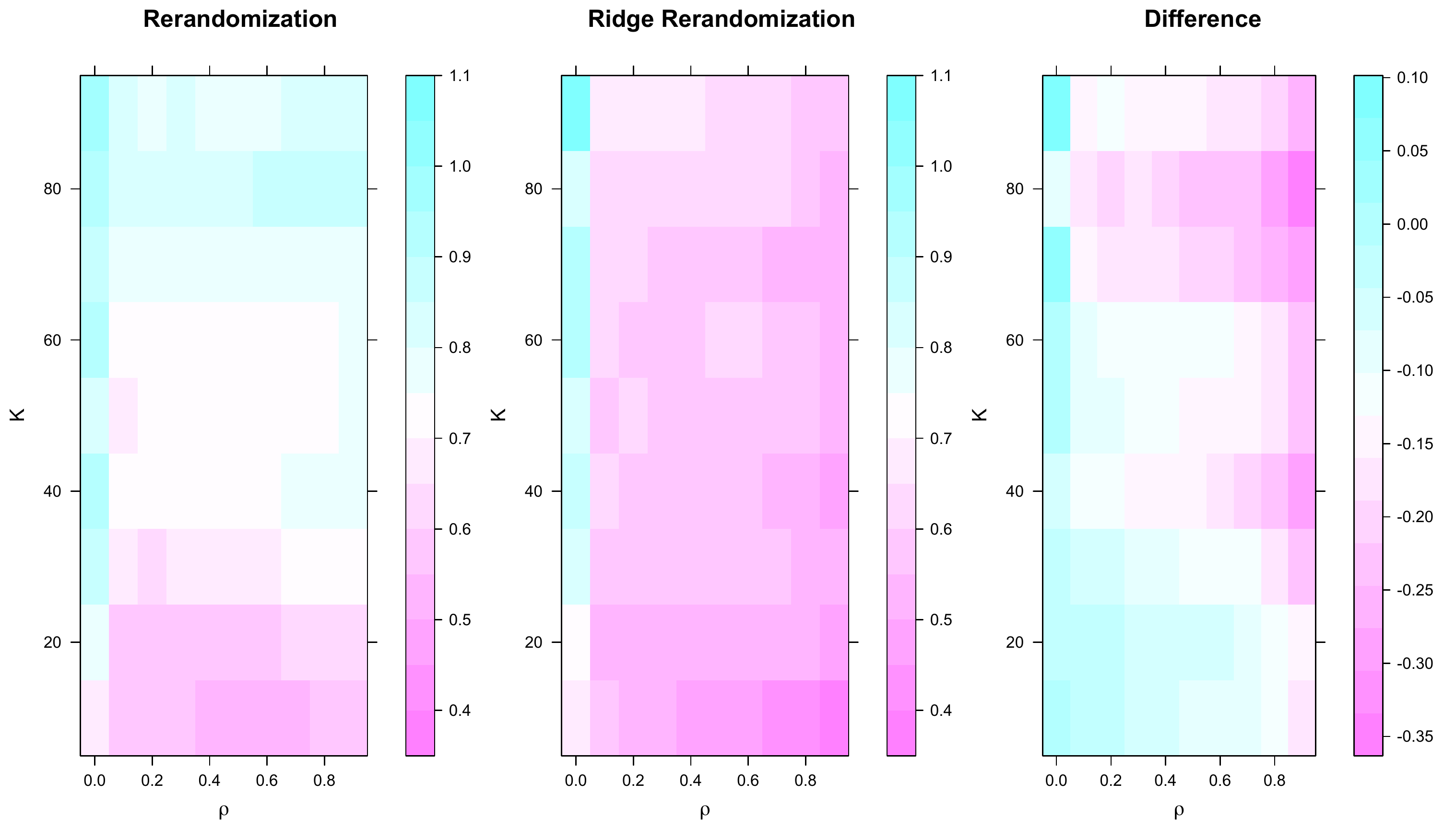}}
	\caption{Relative MSE of $\hat{\tau}$ under rerandomization and ridge rerandomization (relative to randomization) when $\boldsymbol{\beta} = \mathbf{1}_K$ in (\ref{eqn:manyCovariatesPotentialOutcomesMisspecifiedModel}), as well as the difference in relative MSE between the two (i.e., the second plot minus the first).}
	\label{fig:ridgeVSrerandTreatmentEffectMSEMisspecified}
\end{figure}

\begin{figure}[H]
	\centerline{\includegraphics[width=\linewidth]{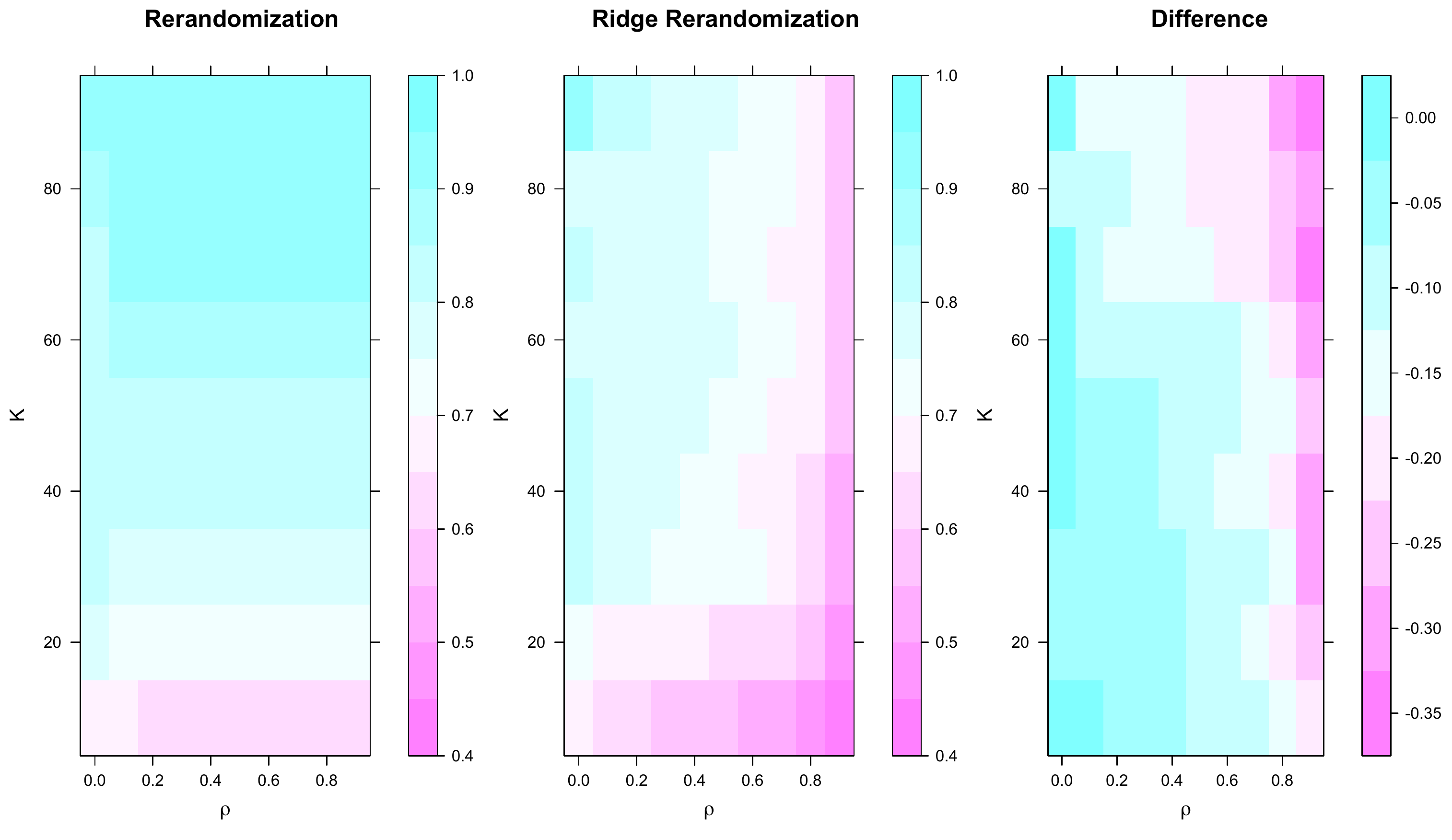}}
	\caption{Relative average confidence interval width under rerandomization and ridge rerandomization (relative to randomization) when $\boldsymbol{\beta} = \mathbf{1}_K$ in (\ref{eqn:manyCovariatesPotentialOutcomesMisspecifiedModel}), as well as the difference between the two (i.e., the second plot minus the first).}
	\label{fig:ridgeVSrerandCIMisspecified}
\end{figure}

\subsection{Summary of Simulation Results}

Importantly, the effectiveness of rerandomization or ridge rerandomization in balancing the covariates does not depend on the covariates' relationship with the outcomes. In other words, the variance reduction results in Figure \ref{fig:ridgeVSrerandPlotVarReduction} do not depend on $\boldsymbol{\beta}$, whereas the treatment effect estimation accuracy results in Figures \ref{fig:ridgeVSrerandTreatmentEffectMSE} and \ref{fig:ridgeVSrerandTreatmentEffectMSEBadBeta} and confidence interval results in Figures \ref{fig:ridgeVSrerandCI} and \ref{fig:ridgeVSrerandCIBadBeta} do. From Figure \ref{fig:ridgeVSrerandPlotVarReduction} we see that ridge rerandomization appears to generally be more effective than rerandomization in balancing covariates in high-dimensional or high-collinearity settings, and from Figures \ref{fig:ridgeVSrerandTreatmentEffectMSE} and Figure \ref{fig:ridgeVSrerandCI} we see that this can result in more precise treatment effect estimators and confidence intervals. These results also hold when the outcome is nonlinearly related with the covariates, as discussed in Section \ref{ss:simsNonlinear}. However, from Section \ref{sss:badBeta}, we see that there are cases where rerandomization can perform better than ridge rerandomization in terms of treatment effect estimation. In particular, if the relationship between the covariates and the outcome is strongly in the direction of the last principal component of the covariate space, rerandomization can perform arbitrarily better than ridge rerandomization, especially when there are only a few number of covariates. In general, the comparison between rerandomization and ridge rerandomization depends on the relationship between the covariates and the outcomes, which is typically not known until after the experiment is conducted.

In summary, these simulations suggest that ridge rerandomization is often preferable over rerandomization by targeting the directions that best explain variation in the covariates rather than the covariates themselves. If the covariates are related to the outcomes (linearly or nonlinearly), ridge rerandomization appears to be an appealing experimental design strategy when there are many covariates and/or highly collinear covariates.

\section{Discussion and Conclusion} \label{s:discussion}

The rerandomization literature has focused on experimental design strategies that utilize the Mahalanobis distance. Starting with \cite{morgan2012rerandomization} and continuing with works such as \cite{morgan2015rerandomization}, \cite{branson2016improving}, \cite{ernst2017sequential}, and \cite{li2018asymptotic}, many theoretical results have been established for rerandomization schemes using the Mahalanobis distance. However, the Mahalanobis distance is known to not perform well in high dimensions or when there are strong collinearities among covariates---settings which the current rerandomization literature has not addressed.

To address experimental design settings where there are many covariates or strong collinearities among covariates, we presented a rerandomization scheme that utilizes a modified Mahalanobis distance. This modified Mahalanobis distance inflates the eigenvalues of the covariance matrix of the covariates, thereby automatically placing a hierarchy of importance among the covariates according to their principal components. Such a quantity has remained largely unexplored in the literature. We established several theoretical properties of this modified Mahalanobis distance, as well as properties of a rerandomization scheme that uses it---an experimental design strategy we call ridge rerandomization. These results establish that ridge rerandomization preserves the unbiasedness of treatment effect estimators and reduces the variance of covariate mean differences. If the covariates are related to the outcomes of the experiment, ridge rerandomization will yield more precise treatment effect estimators than randomization. Furthermore, we conducted a simulation study that suggests that ridge rerandomization is often preferable over rerandomization in high-dimensional or high-collinearity scenarios, which is intuitive given ridge rerandomization's connections to ridge regression.

This modified Mahalanobis distance represents a class of rerandomization criteria, which has connections to principal components and the Euclidean distance. This motivates future work for rerandomization schemes that utilize other criteria. In particular, our theoretical results establish that the benefit of our class of rerandomization schemes over typical rerandomization depends on the covariates' relationship with the outcomes, which usually is not known until after the experiment has been conducted. However, if researchers have prior information about the relationship between the covariates and the outcomes, this information may be useful in selecting rerandomization criteria. An interesting line of future work is further exploring other classes of rerandomization criteria, as well as demonstrating how prior outcome information can be used to select useful rerandomization criteria when designing an experiment.

\bibliographystyle{apa-good}
\bibliography{ridgeRerandomizationBib}

\section{Appendix}

\subsection{Proof of Lemma \ref{lemma:ridgeMDDistribution}} \label{ss:proofLemmaRidgeMDDistribution}

Since $\boldsymbol{\Sigma}>0$, it is invertible and we can write
\begin{equation}
(\mathbf{\Sigma} + \lambda\;I_K)^{-1}=\mathbf{\Sigma}^{-\frac{1}{2}}(I_K + \lambda\;\mathbf{\Sigma}^{-1} )^{-1}\mathbf{\Sigma}^{-\frac{1}{2}}\nonumber
\end{equation}
so that
\begin{equation}
M_\lambda = \widetilde{\mathbf{Z}}^\top (I_K + \lambda\;\mathbf{\Sigma}^{-1} )^{-1}\widetilde{\mathbf{Z}}\nonumber
\end{equation}
where $\widetilde{\mathbf{Z}}=\mathbf{\Sigma}^{-\frac{1}{2}}(\bar{\mathbf{x}}_T - \bar{\mathbf{x}}_C)$. Thanks to the assumed Normality of $\widetilde{\mathbf{Z}}\,|\,\mathbf{x}\sim\mathcal{N}(0,\mathbf{I}_K)$, we may write
\begin{equation*}
M_\lambda\,|\,\mathbf{x} \;\sim\; {\mathbf{Z}}^\top (I_K + \lambda\;\mathbf{\Sigma}^{-1} )^{-1}{\mathbf{Z}}
\end{equation*}
where $\mathbf{Z} = (Z_1 \dots Z_K)^\top\sim \mathcal{N}(0,\mathbf{1}_K)$ marginally and independently of $\mathbf{x}$. The matrix $(I_K + \lambda\;\mathbf{\Sigma}^{-1})^{-1}$ shares the same orthonormal basis x of eigenvectors $\boldsymbol{\Gamma}$ as $\mathbf{\Sigma}$, with corresponding eigenvalues $\lambda_1(\lambda_1+\lambda)^{-1},...,\lambda_K(\lambda_K+\lambda)^{-1}$. As a consequence, we have
\begin{align}
M_\lambda \,|\,\mathbf{x}\;\sim\; (\boldsymbol{\Gamma}^\top\mathbf{Z})^\top\, \textbf{Diag}\left(\left(\frac{\lambda_j}{\lambda_j+\lambda}\right)_{1\leq j \leq K}\right) \,(\boldsymbol{\Gamma}^\top\mathbf{Z})
\label{eq:proof_aux1}
\end{align}
Since $(\boldsymbol{\Gamma}^\top \mathbf{Z})\sim\mathcal{N}(0,\boldsymbol{\Gamma}^\top\boldsymbol{\Gamma})\sim\mathcal{N}(0,\mathbf{I}_K)\sim \mathbf{Z}$ by orthogonality of $\boldsymbol{\Gamma}$, we get
\begin{align*}
	M_\lambda\,|\,\mathbf{x} \;\sim\; \mathbf{Z}^\top \, \textbf{Diag}\left(\left(\frac{\lambda_j}{\lambda_j+\lambda}\right)_{1\leq j \leq K}\right) \,\mathbf{Z} \;=\; \sum_{j=1}^K \frac{\lambda_j}{\lambda_j + \lambda} Z_j^2
\end{align*}
where $Z_1,...,Z_K\stackrel{\text{i.i.d.}}{\sim}\mathcal{N}(0,1)$ and $\lambda_1\geq ...\geq\lambda_K>0$ are the eigenvalues of $\mathbf{\Sigma}$.
\qed

\subsection{Proof of Theorem \ref{thm:ridgeRerandomizationCovariance}} \label{ss:proofTheoremRidgeRerandomizationCovariance}

Using the same notation and reasoning as for the proof of Lemma \ref{lemma:ridgeMDDistribution} in Section \ref{ss:proofLemmaRidgeMDDistribution}, in particular \eqref{eq:proof_aux1}, we can write
\begin{align}
&\cov{\bar{\mathbf{x}}_T - \bar{\mathbf{x}}_C\,|\,\mathbf{x},M_\lambda\leq a_\lambda} \nonumber
\\
=\;\; &\cov{\mathbf{\Sigma}^{1/2}\mathbf{Z}\,\left|\,\mathbf{x},\sum_{j=1}^K \frac{\lambda_j}{\lambda_j + \lambda}(\boldsymbol{\Gamma}^\top\mathbf{Z})^2_j\leq a_\lambda\right.} \nonumber
\\
=\;\; &\cov{\mathbf{\Gamma}\diag{\sqrt{\lambda_{1:K}}}(\boldsymbol{\Gamma}^\top\mathbf{Z})\,\left|\,\sum_{j=1}^K \frac{\lambda_j}{\lambda_j + \lambda}(\boldsymbol{\Gamma}^\top\mathbf{Z})^2_j\leq a_\lambda\right.}
\label{eq:proof42_step1}
\\
=\;\; &\mathbf{\Gamma}\diag{\sqrt{\lambda_{1:K}}}\,\cov{(\boldsymbol{\Gamma}^\top\mathbf{Z})\,\left|\,\sum_{j=1}^K \frac{\lambda_j}{\lambda_j + \lambda}(\boldsymbol{\Gamma}^\top\mathbf{Z})^2_j\leq a_\lambda\right.}\,\diag{\sqrt{\lambda_{1:K}}}\mathbf{\Gamma}^\top
\nonumber
\\
=\;\; &\mathbf{\Gamma}\diag{\sqrt{\lambda_{1:K}}}\,\cov{\mathbf{Z}\,\left|\,\sum_{j=1}^K \frac{\lambda_j}{\lambda_j + \lambda}Z^2_j\leq a_\lambda\right.}\,\diag{\sqrt{\lambda_{1:K}}}\mathbf{\Gamma}^\top
\label{eq:proof42_step2}
\end{align}
where \eqref{eq:proof42_step1} follows from the definition of $\mathbf{\Sigma}^{1/2} = \mathbf{\Gamma}\diag{\sqrt{\lambda_{1:K}}}\mathbf{\Gamma}^\top$ along with the constructed independence of $\mathbf{Z}$ and $\mathbf{x}$ to get rid of the conditioning on $\mathbf{x}$, and \eqref{eq:proof42_step2} follows from $(\boldsymbol{\Gamma}^\top\mathbf{Z})\sim \mathbf{Z}$ by orthogonality of $\boldsymbol{\Gamma}$ and standard Normality of $\mathbf{Z}$. All that is left now is to compute the conditional covariance matrix appearing in \eqref{eq:proof42_step2}. Starting by its diagonal elements, the symmetry of the Normal distribution ensures that $\mathbf{Z}\sim -\mathbf{Z}$, which implies
\begin{align*}
	\mathbb{E}\left[{Z}_k \,\left|\, \sum_{j=1}^K \frac{\lambda_j}{\lambda_j + \lambda} {Z}_j^2 \leq a_{\lambda} \right.\right] &\;=\; \mathbb{E}\left[-{Z}_k \,\left|\, \sum_{j=1}^K \frac{\lambda_j}{\lambda_j + \lambda} (-{Z}_j)^2 \leq a_{\lambda} \right.\right] \\
	&\;=\; - \mathbb{E}\left[{Z}_k \,\left|\, \sum_{j=1}^K \frac{\lambda_j}{\lambda_j + \lambda} ({Z}_j)^2 \leq a_{\lambda} \right.\right]
\end{align*}
for all $k=1,...,K$, so that
\begin{equation*}
	\mathbb{E}\left[{Z}_k \,\left|\, \sum_{j=1}^K \frac{\lambda_j}{\lambda_j + \lambda} {Z}_j^2 \leq a_{\lambda} \right.\right] \;=\; 0
\end{equation*}
Thus, the diagonal elements $d_{k,\lambda}$ of $\text{Cov} \left({\mathbf{Z}} \,\left|\,\sum_{j=1}^K \frac{\lambda_j}{\lambda_j + \lambda} {Z}_j^2 \leq a_{\lambda} \right.\right)$ are given by
\begin{equation}
	d_{k,\lambda} = \text{Var} \left({Z}_k^2 \,\left|\, \sum_{j=1}^K \frac{\lambda_j}{\lambda_j + \lambda}{Z}_j^2 \leq a_{\lambda} \right.\right) = \mathbb{E}\left[ {Z}_k^2 \,\left|\, \sum_{j=1}^K \frac{\lambda_j}{\lambda_j + \lambda} {Z}_j^2 \leq a_{\lambda} \right.\right] 
	\label{eqn:diagonalElements}
\end{equation}
for all $k=1,...,K$. Now for the $(\ell,m)$-element of $\text{Cov} \left({\mathbf{Z}} \,\left|\,\sum_{j=1}^K \frac{\lambda_j}{\lambda_j + \lambda} {Z}_j^2 \leq a_{\lambda} \right.\right)$ with $\ell\neq m$, we use again the symmetry of the Normal distribution by noticing that $\mathbf{Z}\sim \mathbf{Z}^*$, where we define $Z^*_i = Z_i$ for all $i\neq \ell$ and $Z^*_\ell = -Z_\ell$, so that
\begin{align*}
	\text{Cov} \left(Z_\ell,Z_m \,\left|\,\sum_{j=1}^K \frac{\lambda_j}{\lambda_j + \lambda} {Z}_j^2 \leq a_{\lambda} \right.\right) &\;=\; 	\text{Cov} \left(Z_\ell^*,Z_m^* \,\left|\,\sum_{j=1}^K \frac{\lambda_j}{\lambda_j + \lambda} ({Z}_j^*)^2 \leq a_{\lambda} \right.\right)
	\\
	&\;=\; -\,\text{Cov} \left(Z_\ell,Z_m \,\left|\,\sum_{j=1}^K \frac{\lambda_j}{\lambda_j + \lambda} {Z}_j^2 \leq a_{\lambda} \right.\right)
\end{align*}
which leads to 
\begin{equation}
\text{Cov} \left(Z_\ell,Z_m \,\left|\,\sum_{j=1}^K \frac{\lambda_j}{\lambda_j + \lambda} {Z}_j^2 \leq a_{\lambda} \right.\right) \;=\; 0
\label{eq:off_diag}
\end{equation}
for all $1\leq \ell,m\leq K$ such that $\ell\neq m$. Combining \eqref{eqn:diagonalElements} and \eqref{eq:off_diag} gives
\begin{equation}
	\text{Cov} \left({\mathbf{Z}} \,\left|\,\sum_{j=1}^K \frac{\lambda_j}{\lambda_j + \lambda} {Z}_j^2 \leq a_{\lambda} \right.\right) = \textbf{Diag}\left(\left(d_{k,\lambda}\right)_{1\leq k\leq K}\right)
	\label{eq:conditional_cov}
\end{equation}
Plugging\eqref{eq:conditional_cov} back into \eqref{eq:proof42_step2} finally yields
\begin{equation*}
	\cov{\bar{\mathbf{x}}_T - \bar{\mathbf{x}}_C\,|\,\mathbf{x},M_\lambda\leq a_\lambda} = \boldsymbol{\Gamma} {\textbf{Diag}}((\lambda_k\, d_{k,\lambda})_{1\leq k \leq K}) \boldsymbol{\Gamma}^\top.
\end{equation*}
where the $d_{k,\lambda}$'s are given by \eqref{eqn:diagonalElements}. From the expression of $d_{k,\lambda}$, we immediately have $d_{k,\lambda}>0$ for all $k=1,...,K$. By using Equation (13) from \citet{palombi2013note}, we also get 
\begin{equation*}
	\mathbb{E}\left[ {Z}_k^2 \,\left|\, \sum_{j=1}^K \frac{\lambda_j}{\lambda_j + \lambda} {Z}_j^2 \leq a_{\lambda} \right.\right] \;\;<\;\; \mathbb{E}\left[ {Z}_k^2\right] \;=\;1
\end{equation*}
for all $k=1,...,K$. Therefore, we have $d_{k,\lambda}\in(0,1)$ for all $k=1,...,K$.
\qed

\subsection{Proof of Corollary \ref{cor:varReductionRidgeRerandomization}} \label{ss:proofTheoremVarianceReductionBound}
By definition of $v_{k,\lambda}$ and by Theorem \ref{thm:ridgeRerandomizationCovariance}, we have
\begin{align*}
v_{k,\lambda} = \frac{{\text{Var}}\left((\bar{\mathbf{x}}_{T} - \bar{\mathbf{x}}_{C})_k\,|\,\mathbf{x},M_\lambda \leq a_\lambda\right)}{ {\text{Var}}\left((\bar{\mathbf{x}}_{T} - \bar{\mathbf{x}}_{C})_k\,|\,\mathbf{x}\right)} &= \frac{{\text{Cov}}(\bar{\mathbf{x}}_T - \bar{\mathbf{x}}_C\,|\,\mathbf{x},M_{\lambda}\leq a_{\lambda})_{kk}}{{\text{Cov}}(\bar{\mathbf{x}}_T - \bar{\mathbf{x}}_C\,|\,\mathbf{x})_{kk}} 
\\
&= \frac{\left(\boldsymbol{\Gamma} {\textbf{Diag}}\left((\lambda_j \,d_{j,\lambda})_{1\leq j \leq K}\right)\boldsymbol{\Gamma}^\top\right)_{kk}}{\mathbf{\Sigma}_{kk}}.
\end{align*}
Since $\lambda_j(1-d_{j,\lambda}) > 0$ for all $j=1,...,K$, the matrix
\begin{equation*}
	\mathbf{\Sigma} - \boldsymbol{\Gamma} {\textbf{Diag}}\left((\lambda_j \,d_{j,\lambda})_{1\leq j \leq K}\right)\boldsymbol{\Gamma}^\top = \boldsymbol{\Gamma} {\textbf{Diag}}\left((\lambda_j \,(1-d_{j,\lambda}))_{1\leq j \leq K}\right)\boldsymbol{\Gamma}^\top
\end{equation*}
is positive definite. This implies that 
\begin{equation}
	\mathbf{v}^\top \left(\mathbf{\Sigma} - \boldsymbol{\Gamma} {\textbf{Diag}}\left((\lambda_j \,d_{j,\lambda})_{1\leq j \leq K}\right)\boldsymbol{\Gamma}^\top\right) \mathbf{v} \;\;>\;\; 0
	\label{eq:positive_def}
\end{equation}
for all $\mathbf{v}\in\mathbb{R}^K\backslash\{0\}$. In particular, by using \eqref{eq:positive_def} with $\mathbf{v}$ chosen to be the $k$-th canonical basis vector of $\mathbb{R}^K$ (whose elements are all $0$ except its $k$-th element equal to $1$), we get, for all $k=1,...,K$,
\begin{equation}
	\mathbf{\Sigma}_{kk}\;\;>\;\;\left(\boldsymbol{\Gamma} {\textbf{Diag}}\left((\lambda_j \,d_{j,\lambda})_{1\leq j \leq K}\right)\boldsymbol{\Gamma}^\top\right)_{kk}.
\end{equation}
These terms being strictly positive, this leads to $v_{k,\lambda}\in(0,1)$ for all $j=1,...,K$, i.e.
\begin{equation*}
{{\text{Var}}\left((\bar{\mathbf{x}}_{T} - \bar{\mathbf{x}}_{C})_k\,|\,\mathbf{x},M_\lambda \leq a_\lambda\right)}\;\;<\;\;{ {\text{Var}}\left((\bar{\mathbf{x}}_{T} - \bar{\mathbf{x}}_{C})_k\,|\,\mathbf{x}\right)}
\end{equation*}
\qed

\subsection{Proof of Theorem \ref{thm:variancePercentReductionRidgeRerandomizationComparison}} \label{ss:proofTheoremVariancePercentReductionRidgeRerandomizationComparison}
By using \eqref{eq:regression_outcome}, we can write
\begin{equation}
	\hat{\tau} \;=\; (\bar{\mathbf{y}}_T - \bar{\mathbf{y}}_C) \;=\; \tau + \boldsymbol{\beta}^\top (\bar{\mathbf{x}}_T-\bar{\mathbf{x}}_C) + (\bar{\boldsymbol{\epsilon}}_T-\bar{\boldsymbol{\epsilon}}_C)
\end{equation}
By conditional independence of $(\bar{\boldsymbol{\epsilon}}_T-\bar{\boldsymbol{\epsilon}}_C)$ and $(\bar{\mathbf{x}}_T-\bar{\mathbf{x}}_C)$ given $\mathbf{x}$, we have
\begin{align}
\text{Var}(\hat{\tau}\,|\,\mathbf{x}) \;&=\; \text{Var}(\boldsymbol{\beta}^\top (\bar{\mathbf{x}}_T-\bar{\mathbf{x}}_C)\,|\,\mathbf{x}) + \text{Var}(\bar{\boldsymbol{\epsilon}}_T-\bar{\boldsymbol{\epsilon}}_C\,|\,\mathbf{x}) \nonumber
\\
\;&=\; \boldsymbol{\beta}^\top\boldsymbol{\Sigma}\boldsymbol{\beta} + \text{Var}(\bar{\boldsymbol{\epsilon}}_T-\bar{\boldsymbol{\epsilon}}_C\,|\,\mathbf{x}) 
\label{eq:var_noise}
\end{align}
Conditional on $\mathbf{x}$, $M_\lambda$ is a deterministic function of $(\bar{\mathbf{x}}_T-\bar{\mathbf{x}}_C)$, thus $(\bar{\boldsymbol{\epsilon}}_T-\bar{\boldsymbol{\epsilon}}_C)$ is conditionally independent of $M_\lambda$ given $\mathbf{x}$. This leads to 
\begin{align}
\text{Var}(\hat{\tau}\,|\,\mathbf{x},M_\lambda\leq a_\lambda) \;&=\; \text{Var}(\boldsymbol{\beta}^\top (\bar{\mathbf{x}}_T-\bar{\mathbf{x}}_C)\,|\,\mathbf{x},M_\lambda\leq a_\lambda) + \text{Var}(\bar{\boldsymbol{\epsilon}}_T-\bar{\boldsymbol{\epsilon}}_C\,|\,\mathbf{x},M_\lambda\leq a_\lambda) \nonumber
\\
\;&=\; \boldsymbol{\beta}^\top\text{Cov}(\bar{\mathbf{x}}_T-\bar{\mathbf{x}}_C\,|\,\mathbf{x},M_\lambda\leq a_\lambda)\boldsymbol{\beta} + \text{Var}(\bar{\boldsymbol{\epsilon}}_T-\bar{\boldsymbol{\epsilon}}_C\,|\,\mathbf{x}) 
\label{eq:step_cond_indep}
\\
\;&=\;\boldsymbol{\beta}^\top \boldsymbol{\Gamma} {\textbf{Diag}} \left((\lambda_k d_{k,\lambda})_{1\leq k \leq K} \right) \boldsymbol{\Gamma}^\top \boldsymbol{\beta} + \text{Var}(\bar{\boldsymbol{\epsilon}}_T-\bar{\boldsymbol{\epsilon}}_C\,|\,\mathbf{x}) 
\label{eq:step_thm}
\end{align}
where \eqref{eq:step_cond_indep} follows from the conditional independence of $(\bar{\boldsymbol{\epsilon}}_T-\bar{\boldsymbol{\epsilon}}_C)$ and $M_\lambda$ given $\mathbf{x}$, and \eqref{eq:step_thm} follows from Theorem \ref{thm:ridgeRerandomizationCovariance}. By plugging \eqref{eq:var_noise} into \eqref{eq:step_thm}, we get
\begin{align*}
{\text{Var}}(\hat{\tau}\,|\,\mathbf{x}) - {\text{Var}}(\hat{\tau}\,|\,\mathbf{x},M_{\lambda} \leq a_{\lambda}) \;&=\; \boldsymbol{\beta}^\top (\boldsymbol{\Sigma}-\boldsymbol{\Gamma}{\textbf{Diag}} \left((\lambda_k d_{k,\lambda} )_{1\leq k \leq K} \right) \boldsymbol{\Gamma}^\top) \boldsymbol{\beta}
\\
\;&=\;\boldsymbol{\beta}^\top \boldsymbol{\Gamma} {\textbf{Diag}} \left((\lambda_k \left(1 - d_{k,\lambda} \right))_{1\leq k \leq K} \right) \boldsymbol{\Gamma}^\top \boldsymbol{\beta}
\end{align*}
As explained by \eqref{eq:positive_def} in the proof of Corollary \ref{cor:varReductionRidgeRerandomization}, the positive definiteness of the matrix $\boldsymbol{\Gamma}{\textbf{Diag}} \left((\lambda_k \left(1 - d_{k,\lambda} \right))_{1\leq k \leq K} \right) \boldsymbol{\Gamma}^\top$ guarantees that
\begin{align*}
 {\text{Var}}(\hat{\tau}\,|\,\mathbf{x},M_{\lambda} \leq a_{\lambda})\;\;\leq\;\;{\text{Var}}(\hat{\tau}\,|\,\mathbf{x})
\end{align*}
for all $\boldsymbol{\beta}\in\mathbb{R}^K$, with equality if and only if  $\boldsymbol{\beta}=0$.
\qed

\end{document}